\documentclass[a4paper,reqno,10pt]{amsart}

\usepackage[USenglish]{babel}
\usepackage[utf8x]{inputenc}
\usepackage[style=american]{csquotes}
\usepackage{ucs}
\usepackage{yfonts}
\usepackage{amsmath,amsthm,amscd,amssymb}
\usepackage{MnSymbol}
\usepackage[all]{xy}
\usepackage{smartref}
\usepackage{algorithmic}
\usepackage{tikz}
\usepackage{algorithm, caption}
\usepackage{hyperref}
\usepackage[lite,shortalphabetic]{amsrefs}

\theoremstyle{definition}
\newtheorem{defn}{Definition}[section]
\newtheorem{ex}[defn]{Example}

\newtheorem{remark}[defn]{Remark}

\newtheorem*{acknowledgement}{Acknowledgement}

\theoremstyle{plain}
\newtheorem{theorem}[defn]{Theorem}
\newtheorem{corollary}[defn]{Corollary}
\newtheorem{lemma}[defn]{Lemma}
\newtheorem{prop}[defn]{Proposition}

\newcommand{\abs}[1]{\left\lvert #1 \right\rvert} 
\newcommand{\gnrt}[1]{\left\langle #1 \right\rangle} 

\newcommand{\Z}{\mathbb{Z}}

\newcommand{\N}{\mathbb{N}}
\newcommand{\R}{\mathbb{R}}

\newcommand{\symm}{\mathbb{S}}

\renewcommand{\Im}{\textnormal{Im}}
\newcommand{\relint}{\textnormal{rel int}\,}

\newcommand{\Star}{\textnormal{Star}}

\newcommand{\conv}{\textnormal{conv}}

\newcommand{\mn}{\mathcal{M}_{n}}

\newcommand{\val}{\textnormal{val}}
\newcommand{\ft}{\textnormal{ft}}
\newcommand{\curly}[1]{\mathcal{#1}}

\newcommand{\polymake}{\texttt{polymake}}
\newcommand{\atint}{\texttt{a-tint}}

\newcommand{\tsr}{\otimes} 
\newcommand{\wo}{\setminus} 

\definecolor{comment}{rgb}{.2,.2,.2}
\newcommand{\INPUT}{\STATE \textbf{Input: }}
\newcommand{\OUTPUT}{\STATE \textbf{Output: }}
\newcommand{\COMM}[1]{\STATE //{\it \textcolor{comment}{#1}}}

\newcommand{\atcommand}{atint  $>$ }
\renewcommand{\>}{$>$}

\newcommand{\newl}{\indent\par}

\newcommand{\emphex}[1]{\emph{#1}\index{#1}} 

\numberwithin{equation}{section}
\newcommand{\arxiv}[1]{%
  \href{http://arxiv.org/abs/#1}{arxiv:#1}%
}

\begin{document}
\title[a-tint: algorithmic tropical intersection theory]{a-tint: a polymake extension for algorithmic tropical intersection theory}
\author {Simon Hampe}
\address {Simon Hampe, Fachbereich Mathematik, Technische Universit\"{a}t Kaiserslautern, Postfach 3049, 67653 Kaiserslautern, Germany}
\email {hampe@mathematik.uni-kl.de}
\begin{abstract}
  In this paper we study algorithmic aspects of tropical intersection theory. We analyse how divisors and intersection products on tropical cycles can actually be computed using polyhedral geometry. The main focus of this paper is the study of moduli spaces, where the underlying combinatorics of the varieties involved allow a much more efficient way of computing certain tropical cycles. The algorithms discussed here have been implemented in an extension for \polymake, a software for polyhedral computations.
\end{abstract}
\maketitle

\section{Introduction}

Tropical intersection theory has proved to be a powerful tool in tropical geometry. The basic ideas for an intersection theory in $\R^n$ based on divisors of rational functions were first laid out in \cite{mikh} and were then further developed in \cite{AR}. Even earlier, \cite{fulton} had proved the \emph{fan displacement rule} in the context of toric geometry. It describes how cohomology classes of a toric variety $X(\Delta)$ can be multiplied using a generic translation of $\Delta$. This was later translated to the concept of \emph{stable intersections} of tropical varieties. One can force two tropical varieties to intersect in the correct dimension by translating one of them locally by a generic vector. The intersection multiplicities are then computed using the formula from the fan displacement rule. 

An intersection product in matroid fans was introduced in \cite{shaw},\cite{francoisrau}. In particular, this made it now possible to do intersection theory on moduli spaces. 

Tropical intersection theory has many applications. For example, one can use \cite{fulton} to see that certain intersection products in toric varieties can be computed as tropical intersection products. It has also been used in \cite{obstructions} to study the relative realizability of tropical curves. A prominent example of the usefulness of tropical intersection theory is enumerative geometry (see for example \cite{GKM}, \cite{raumoduli}, \cite{psiclasses}). One can formulate many combinatorially complex enumerative problems in terms of much simpler intersection products on moduli spaces. 

However, in concrete cases these products are still tedious to compute by hand, especially in higher dimensions. For the purpose of testing new conjectures or studying examples of a new and unfamiliar object, one is often interested in such explicit computations. This paper aims to analyse how one can efficiently compute divisors, products of tropical cycles and other constructions frequently occurring in tropical intersection theory.

After briefly discussing the basic notions of polyhedral and tropical geometry, we study some basic operations in tropical geometry in Section \ref{section_basic}: We describe how a lattice normal vector and the divisor of a rational function are computed and we give an algorithm that can check whether a given tropical cycle is irreducible. We prove that this property is independent of the actual polyhedral structure of the cycle, then show how it can be computed for a given structure.

In Section \ref{section_intersection} we review the possible definitions of a tropical intersection product and discuss how these definitions give rise to different algorithms to compute an intersection product. The definition by Allermann and Rau (\cite{AR}) immediately suggests an algorithm, but the definition from \cite{jensenyu}, which is closer to the fan displacement rule (\cite{fulton}), provides a much more efficient method.

In Section \ref{section_matroid} we give two different, but equivalent definitions of a \emph{Bergman fan} associated to a matroid. The first describes this fan as a subfan of the normal fan of the matroid polytope, which gives an immediate, but very inefficient method to compute it. The second definition in terms of circuits has been used by Felipe Rincón in his software \texttt{TropLi}(\cites{tropli,troplisoftware}) to compute Bergman fans of matrix matroids. \atint\ uses its own implementation of Rincón's algorithm, which is slower than the original program, but still much faster than the normal fan algorithm.

Section \ref{section_moduli} is the main focus of this paper: Here we discuss how the combinatorial structure of the moduli spaces $\mn$ of rational $n$-marked tropical curves can be used to efficiently compute these spaces or certain subcycles thereof. We introduce Prüfer sequences and show how they are in bijection to combinatorial types of rational curves, using a slight variation of Prüfer's original argument for trees in complete graphs. Since Prüfer sequences are relatively easy to enumerate, we can make use of this to compute (parts of) $\mn$. We use this method to construct an algorithm computing products of tropical Psi-classes. We also show how the combinatorics of a rational curve can be retrieved from its metric vector by an algorithm of Buneman (\cite{buneman}) and we study the local structure of $\mn$. More precisely, we show that locally around each point $\mn$ is the Cartesian product of some $\R^k$ and several $\mathcal{M}_i$ with $i \leq n$.

All the algorithms described in this paper have been implemented by the author in \atint, an extension for \polymake. The latter is a software for polyhedral computations (\cite{polymake}). More information, as well as packages and source code for installation can be found at \url{www.polymake.org}. 
In Section \ref{section_appendix} we list some open questions and the main features of \atint. We also include benchmarking tables that show how some algorithms compare to one another or react to a change in certain parameters.

You can obtain \atint\ under \url{https://bitbucket.org/hampe/atint}, installation instructions and a user manual are provided under the \texttt{Wiki} link. All the \textbf{polymake examples} in this paper can easily be reproduced: Having installed \polymake\ and \atint, enter \texttt{polymake -Aatint} in your console, then type the commands from the example.

\begin{acknowledgement}
The author is supported by the Deutsche Forschungsgemeinschaft grants GA 636 / 4-1, MA 4797/3-1. The software project \atint\ is part of the DFG priority project SPP 1489 (\url{www.computeralgebra.de}).

I would like to thank Andreas Gathmann, Anders Jensen, Hannah Markwig, Thomas Markwig and Kirsten Schmitz for their support and many helpful discussions.
\end{acknowledgement}

\section{Preliminaries: Polyhedral and tropical geometry}\label{section_intro}

In this section we establish the basic terms and definitions of polyhedral and tropical geometry needed in this paper. For a more thorough introduction to polytopes and polyhedra see for example \cite{ziegler}.

\subsection{Polyhedra and polyhedral complexes}

A \emph{polyhedron} or \emph{polyhedral cell} in $V = \R^n$ is a set of the form $$\sigma = \{x : Ax \geq b\}$$ where $A \in \R^{m \times n}, b \in \R^m$, i.e.\ it is an intersection of finitely many halfspaces. We call $\sigma$ a \emph{cone} if $b = 0$. 

Equivalently, any polyhedron $\sigma$ can be described as $$\sigma = \conv\{p_1,\dots,p_k\} + \R_{\geq 0}r_1 + \dots + \R_{\geq 0}r_l + L$$
where $p_1,\dots,p_k,r_1,\dots,r_l \in \R^n$, $L$ is a linear subspace of $\R^n$ and $+$ denotes the Minkowski sum of sets:
$$A + B = \{a+b; a \in A, b \in B\}.$$ The first description is often called an $\curly{H}$-description of $\sigma$ and the second is a $\curly{V}$-description of $\sigma$. 

It is a well known algorithmic problem in convex geometry to compute one of these descriptions from the other. In fact, both directions are computationally equivalent and there are several well-known \emph{convex-hull-algorithms}. Most notable are the double-description method \cite{doubledesc}, the reverse search method \cite{lrs} and the beneath-and-beyond algorithm (e.g.\ \cite{gruenb},\cite{edelsbr}). Generally speaking, each of these algorithms behaves very well in terms of complexity for a certain class of polyhedra, but very badly for some other types (see \cite{ABS} for a more in-depth discussion of this). Since in tropical geometry all kinds of polyhedra can occur, it is very difficult to pick an optimal algorithm. It is also an open problem, whether there exists a convex hull algorithm with polynomial complexity in both input and output. All of the above mentioned algorithms are implemented in \polymake \ or in  libraries used by \polymake. At the moment, all algorithms in \atint\ use the implementation of the double-description algorithm by Fukuda \cite{fukuda_cdd}. 

For any polyhedron $\sigma$ we denote by $V_\sigma$ the vector space associated to the affine space spanned by $\sigma$, i.e.\ $$V_\sigma := \gnrt{a-b; a,b \in \sigma}$$ We denote by $\Lambda_\sigma := V_\sigma \cap \Z^n$ its associated lattice. The \emph{dimension} of $\sigma$ is the dimension of $V_\sigma$.

A \emph{face} of $\sigma$ is any subset $\tau$ that can be written as $\sigma \cap H$, where $H = \{x: \gnrt{x,a} = \lambda\}$ is an affine hyperplane such that $\sigma$ is contained in one of the halfspaces $\{x: \gnrt{x,a} \geq \lambda\}$ or $\{x: \gnrt{x,a} \leq \lambda\}$ (i.e.\ we change one or more of the inequalities defining $\sigma$ to an equality). If $\tau$ is a face of $\sigma$, we write this $\tau \leq \sigma$ (or $\tau < \sigma$ if the inclusion is proper). By convention we will also say that $\sigma$ is a face of itself.

Finally, the \emph{relative interior} of a polyhedron is the set
$$\relint(\sigma) := \sigma \wo \bigcup_{\tau < \sigma} \tau$$
A \emph{polyhedral complex} is a set $\Sigma$ of polyhedra that fulfills the following properties:
\begin{itemize}
 \item For each $\sigma \in \Sigma$ and each face $\tau \leq \sigma$, $\tau \in \Sigma$
 \item For each two $\sigma, \sigma' \in \Sigma$, the intersection is a face of both.
\end{itemize}
If all of the polyhedra in $\Sigma$ are cones, we call $\Sigma$ a \emph{fan}.

We will denote by $\Sigma^{(k)}$ the set of all $k$-dimensional polyhedra in $\Sigma$ and set the dimension of $\Sigma$ to be the largest dimension of any polyhedron in $\Sigma$. The set-theoretic union of all cells in $\Sigma$ is denoted by $\abs{\Sigma}$, the \emph{support} of $\Sigma$. We call $\Sigma$ \emph{pure-dimensional} or \emph{pure} if all inclusion-maximal cells are of the same dimension. We call $\Sigma$ \emph{rational} if all polyhedral cells are defined by inequalities $Ax \geq b$ with rational coefficients $A$. If not explicitly stated otherwise, all complexes and fans in this paper will be pure and rational.

Note that a polyhedral complex is uniquely defined by giving all its top-dimensional cells. Hence we will often identify a polyhedral complex with its set of maximal cells.

The \emph{cartesian product} of two polyhedral complexes $\Sigma$ and $\Sigma'$ is the polyhedral complex
$$\Sigma \times \Sigma' := \{\sigma \times \sigma'; \sigma \in \Sigma, \sigma' \in \Sigma'\}.$$

The last definition we need is the \emph{normal fan} of a \emph{polytope}, i.e.\ a bounded polyhedron: Let $\sigma$ be a polytope, $\tau$ any face of $\sigma$. The \emph{normal cone} of $\tau$ in $\sigma$ is  
$$N_{\tau,\sigma} := \overline{\{w \in \R^n: \gnrt{w,t} = \max\{\gnrt{w,x}; x \in \sigma\} \textnormal{ for all } t \in \tau\}}$$
i.e.\ the closure of the set of all linear forms which take their maximum on $\tau$. These sets are in fact cones forming a fan and the collection of these cones is called the \emph{normal fan} $N_\sigma$ of $\sigma$.

\subsection{Tropical geometry}

Let $X$ be a pure $d$-dimensional rational polyhedral complex in $\R^n$. Let $\sigma \in X^{(d)}$ and assume $\tau \leq \sigma$ is a face of dimension $d-1$. The \emph{primitive normal vector} of $\tau$ with respect to $\sigma$ is defined as follows: By definition there is a linear form $g \in (\Z^n)^\vee$ such that its minimal locus on $\sigma$ is $\tau$. Then there is a unique generator of $\Lambda_\sigma/\Lambda_\tau \cong \Z$, denoted by $u_{\sigma/\tau}$, such that $g(u_{\sigma/\tau}) > 0$. One can also define this for a polyhedral complex in a vector space $\Lambda \tsr \R$, for a prescribed lattice $\Lambda$ (see for example \cite{GKM}). If not stated otherwise, we will however always consider the standard lattice $\Z^n$.

A \emph{tropical cycle} $(X,\omega)$ is a pure rational $d$-dimensional complex $X$ together with a \emph{weight function} $\omega: X^{(d)} \to \Z$ such that for all codimension one faces $\tau \in X^{(d-1)}$ it fulfills the \emph{balancing condition}:
$$\sum_{\sigma > \tau} \omega(\sigma) u_{\sigma/\tau} = 0 \in V / V_\tau$$
We call $X$ a \emph{tropical variety} if furthermore all weights are positive. 

The \emph{cartesian product} of two tropical cycles $X,Y$ is the cartesian product of their underlying polyhedral complexes, equipped with the weight function
$$\omega_{X \times Y}: \sigma \times \sigma' \mapsto \omega_X(\sigma) \cdot \omega_Y(\sigma').$$
It is easy to see that this is again a tropical cycle.

Two tropical cycles are considered \emph{equivalent} if they have the same support and there is a finer polyhedral structure on this support that respects both weight functions. Hence we will sometimes distinguish between a tropical cycle $X$ and a specific polyhedral structure $\curly{X}$.

\begin{figure}
 \centering
 \begin{tikzpicture}
  \matrix[column sep = 10pt]{
    \pgftext{\includegraphics[scale = 0.15]{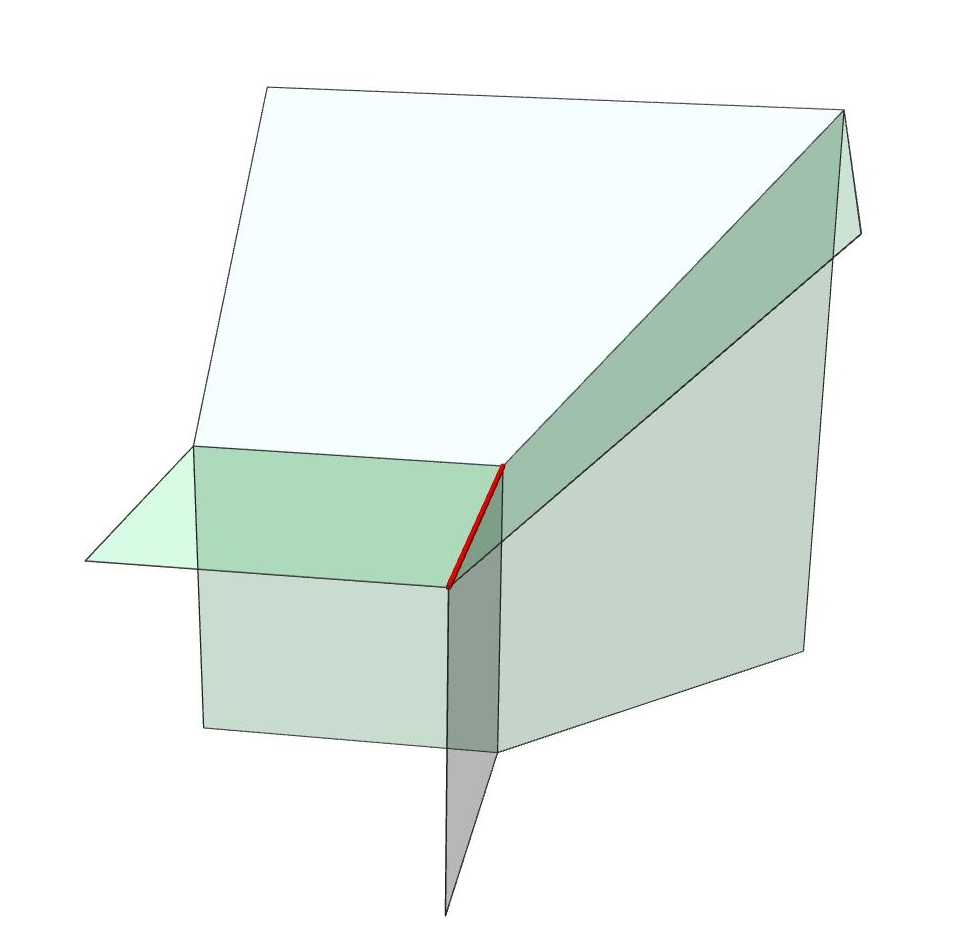}} &
    \draw (0,0) -- (1,1);
    \draw (0,0) -- (-1,0);
    \draw (0,0) -- (0,-1);
    \fill[red] (0,0) circle (2pt);
    \\
  };
 \end{tikzpicture}
  \caption{A tropical plane $L$ and its local picture $\Star_L(\tau)$, where $\tau$ is the codimension one face marked in red.}\label{intro_fig_star}
\end{figure}
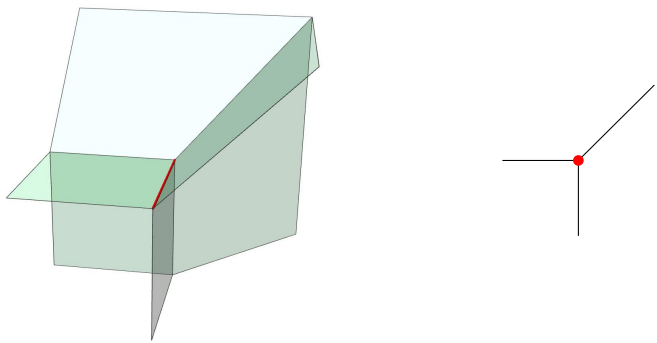

We also want to define the local picture of $\curly{X}$ around a given cell: Let $\tau \in \curly{X}$ be any polyhedral cell. Let $\Pi: V \to V/V_\tau$ be the residue morphism. We define 
$$\Star_\curly{X}(\tau) := \{ \R_{\geq 0} \cdot \Pi(\sigma - \tau); \tau < \sigma \in X\} \cup \{0\}$$
which is a fan in $V/V_\tau$ (on the lattice $\Lambda / \Lambda_\tau$). If we furthermore equip $\Star_\curly{X}(\tau)$ with the weight function $\omega_\Star(\R_{\geq 0} \cdot \Pi(\sigma - \tau)) = \omega_X(\sigma)$ for all maximal $\sigma$, then $(\Star_\curly{X}(\tau), \omega_\Star)$ is a tropical fan cycle. It is now easy to see that a weighted complex $\curly{X}$ is balanced around a codimension cell $\tau$, if and only if the one-dimensional fan $\Star_{\curly{X}}(\tau)$ is balanced. An example for this construction can be found in Figure \ref{intro_fig_star}.

\section{Basic computations in tropical geometry}\label{section_basic}

\subsection{Computing the primitive normal vector} The primitive normal vector $u_{\sigma/\tau}$ defined in the previous section is an essential part of most formulas and calculations in tropical geometry. Hence we will need an algorithm to compute it. An important tool in this computation is the \emph{Hermite normal form} of an integer matrix:

\begin{defn}
 Let $M \in \Z^{m\times n}$ be a matrix with $n \geq m$ and assume $M$ has full rank. We say that $M$ is in \emph{Hermite normal form} (HNF) if it is of the form 
$$M = ( 0_{m \times (n-m)}, T )$$
where $T = (t_{i,j})$ is an upper triangular matrix with $t_{i,i} > 0$ and for $j > i$ we have $t_{i,i} > t_{i,j} \geq 0$.
\end{defn}

\begin{remark}
 We are actually only interested in the fact that $T$ is an upper triangular invertible matrix. Furthermore, it is known that for any $A \in \Z^{m \times n}$ of full rank there exists a $U \in \textnormal{GL}_n(\Z)$ such that $B = AU$ is in HNF (see for example \cite{cohen}*{2.4}).
\end{remark}

\begin{prop}
 Let $X \subset \R^n$ be a $d$-dimensional tropical cycle, $\tau \in X^{(d-1)}$. Let $A \in \Z^{(n-d+1) \times n}$ such that $V_\tau = \ker A$ and $V_\sigma = \ker \tilde{A}$, where $\tilde{A}$ denotes $A$ without its first row. Let $U \in \textnormal{GL}_n(\Z)$ such that
$$AU = (0_{(n-d+1) \times (d-1)}, T)$$
is in HNF. Denote by $U_{*i}$ the $i$-th column of $U$. Then:
\begin{enumerate}
 \item $U_{*1},\dots,U_{*d-1}$ is a lattice basis for $\Lambda_\tau$.
 \item $U_{*1},\dots,U_{*d}$ is a lattice basis for $\Lambda_\sigma$.
\end{enumerate}
In particular $U_{*d} = \pm u_{\sigma/\tau}$ mod $V_\tau$.
\begin{proof}
 It is clear that $U_{*1},\dots, U_{*d-1}$ form an $\R$-basis for $\ker A$ and the fact that $\det U = \pm 1$ ensures that it is a lattice basis. Removing the first row of $A$ corresponds to removing the first row of $AU$ so we obtain an additional column of zeros. Hence $U_{*1},\dots,U_{*d}$ is a basis of $\ker \tilde{A}$ and $U_{*d}$ is a generator of $\Lambda_\sigma / \Lambda_\tau$.
\end{proof}
\end{prop}

\begin{remark}
 In \cite{cohen}*{2.4.3}, Cohen suggests an algorithm for computing the HNF of a matrix using integer Gaussian elimination. However, he already states that this algorithm is useless for practical applications, since the coefficients in intermediate steps of the computation explode too quickly. A more practical solution is an LLL-based normal form algorithm that reduces the coefficients in between elimination steps. \atint\ uses an implementation based on the algorithm designed by Havas, Majevski and Matthews in \cite{havas}.

Note that, knowing the primitive normal vector up to sign, it is easy to determine its final form,  since we know that the linear form defined by $u_{\sigma/\tau}$ must be positive on $\sigma$. So we only have to compute the scalar product of $U_{*d}$ with any ray in $\sigma$ that is not in $\tau$ and check if it is positive.

What remains is to compute the matrix $A$ such that $A = V_\tau$ and $\tilde{A} = V_\sigma$. There is an obvious notion of an \emph{irredundant} $\mathcal{H}$-description of $\tau$: Assume 
$$\tau = \bigcap_{i=1}^r H_i \cap \bigcap_{j=1}^s S_j$$
where $H_i = \{x \in \R^n;\, \gnrt{x,z_i} = \alpha_i\}$ and $S_j = \{x \in \R^n;\, \gnrt{x, w_j} \geq \beta_j\}$ for some $z_i, w_j \in \Z^n, \alpha_i, \beta_j \in \R$. This is considered an \emph{irredundant} $\mathcal{H}$-description if we cannot remove any of these without changing the intersection and we cannot change any of the inequalities into an equality. Note that most convex hull algorithms return such an irredundant description. Now it is basic linear algebra to see the following:
\end{remark}

\begin{lemma}
 Let $\sigma \in \R^n$ be a polyhedron given by an irredundant $\mathcal{H}$-representation $\sigma = (\bigcap_{i=1}^r H_i) \cap (\bigcap_{j=1}^s S_j)$ with $H_i = \{x \in \R^n;\, \gnrt{x ,z_i} = \alpha_i\}$. Denote by $H_i^0 =$ $\left\{x \in \R^n;\, \gnrt{x ,z_i }= 0\right\}$. Then 
$$V_\sigma = \bigcap_{i=1}^r H_i^0 = \ker \begin{pmatrix} z_1 \\ \vdots \\ z_r \end{pmatrix} $$
Furthermore, if $\tau$ is a codimension one face of $\sigma$ with irredundant $\mathcal{H}$-representation $\tau = (\bigcap_{i=1}^{r'} H_i') \cap (\bigcap_{j=1}^{s'} S_j')$, then there is an $l \in \{1,\dots,r'\}$, such that
$$V_\tau = (\bigcap_{i = 1}^r H_i^0) \cap H_l'.$$
\end{lemma}

\begin{remark}
 There is an additional trick that can make lattice normal computations speed up by a considerable factor, especially if the codimension of the tropical variety is large. Assume  you want to compute normal vectors of a 10-dimensional tropical variety in $\R^{20}$. In this case we would have to compute the HNF of $11 \times 20$-matrices. However, for computing $u_{\sigma/\tau}$ we can project $\sigma$ onto $V_\sigma$. Now the matrix of the codimension one face $\tau$ is only a $1 \times 10$-matrix. The normal form of this matrix can of course be computed much faster. Note that we have to take care that the projection induces a lattice isomorphism on $\sigma$. For this, we have to compute a lattice basis of $\sigma$, which still requires computation of an HNF of the matrix associated to $\sigma$ - \emph{but only once}, instead of once for each codimension one face. 
\end{remark}

\begin{algorithm}[ht]
  \caption*{\textbf{polymake example: Computing a tropical variety}.\\ This creates the weighted complex consisting of the four orthants of $\R^2$ with weight 1 and checks if it is balanced. The maximal cones are represented in terms of indices of the rays (starting the count with 0).}
  \flushleft
  
  \texttt{
  \begin{tabular}{l l}
   \atcommand & \$w = new WeightedComplex(\\
	      &RAYS=\>[[1,0],[-1,0],[0,1],[0,-1]],\\
	      &MAXIMAL\_CONES=\>[[0,2],[0,3],[1,2],[1,3]],\\
	      &TROPICAL\_WEIGHTS=\>[1,1,1,1]);\\
   \atcommand & print \$w-\>IS\_BALANCED;\\
   1 
  \end{tabular}}
\end{algorithm}
%
%

\subsection{Divisors of rational functions}

The most basic operation in tropical intersection theory is the computation of the divisor of a rational function. Let us first discuss how we define a rational function and its divisor. Our definition is the same as in \cite{AR}:
\begin{defn}
 Let $X$ be a tropical variety. A \emph{rational function} on $X$ is a function $\varphi: X \to \R$ that is affine linear with integer slope on each cell of some arbitrary polyhedral structure $\mathcal{X}$ of $X$.

The \emph{divisor} of $\varphi$ on $X$, denoted by $\varphi \cdot X$, is defined as follows: Choose a polyhedral structure $\curly{X}$ of $X$ such that $\varphi$ is affine linear on each cell. Let $\curly{X}' = \curly{X}^{(\dim X - 1)}$ be the codimension one skeleton. For each $\tau \in \curly{X}'$, we define its weight via
$$\omega_{\varphi \cdot X}(\tau) = \left(\sum_{\sigma > \tau} \omega(\sigma) \varphi_\sigma(u_{\sigma/\tau})\right) - \varphi_\tau\left(\sum_{\sigma > \tau} \omega(\sigma) u_{\sigma/\tau}\right)$$
where $\varphi_\sigma$ and $\varphi_\tau$ denote the linear part of the restriction of $\varphi$ to the respective cell. Then
$$\varphi \cdot X := (\curly{X}', \omega_{\varphi \cdot X})$$
\end{defn}

\begin{remark}
 While the computation of the weights on the divisor is relatively easy to implement, the main problem is computing the appropriate polyhedral structure. The most general form of a rational function $\varphi$ on some cycle $X$ would be given by its domain, a polyhedral complex $Y$ with $\abs{X} \subseteq \abs{Y}$ together with the values and slopes of $\varphi$ on the vertices and rays of $Y$. To make sure that $\varphi$ is affine linear on each cell of $X$, we then have to compute the intersection of the complexes, which boils down to computing the pairwise intersection of all maximal cones of $X$ and $Y$. Here lies the main problem of computing divisors: One usually computes the intersection of two cones by converting them to an $\curly{H}$-description and converting the joint description back to a $\curly{V}$-description via some convex hull algorithm. But as we discussed earlier, so far no convex hull algorithm is known that has polynomial runtime for all polyhedra. Also, \cite{tiwary} shows that computing the intersection of two $\curly{V}$-polyhedra is NP-complete.

Hence we already see a crucial factor for computing divisors (besides the obvious ones: dimension and ambient dimension): The number of maximal cones of the tropical cycle and the domain of the rational function. Table \ref{table_divisor_poly} in the appendix shows how divisor computation is affected by these parameters.
\end{remark}

\begin{figure}[ht]
 \includegraphics[scale=0.2]{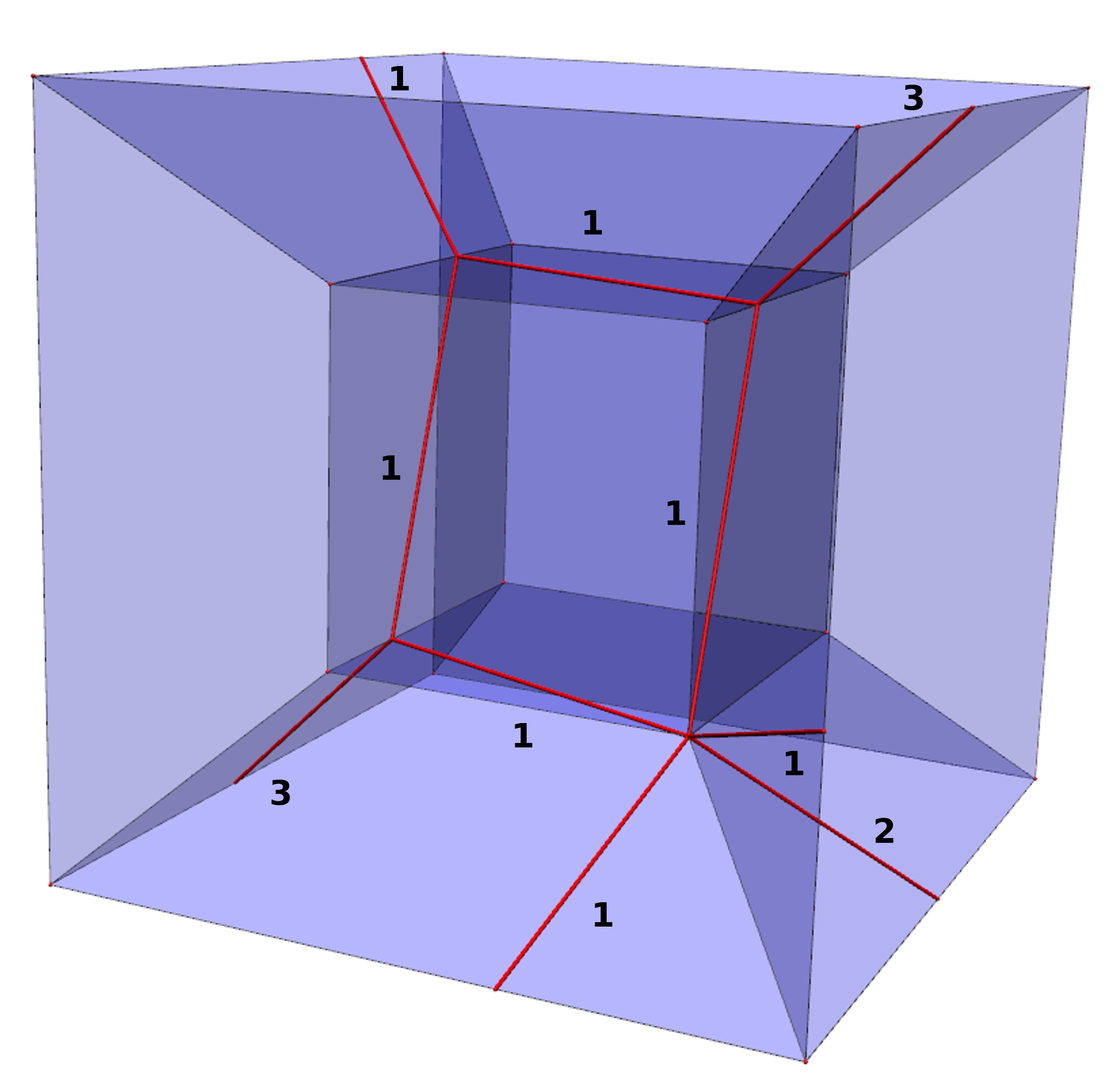}
\caption{The surface is $X = \max\{1,x,y,z,-x,-y,-z\} \cdot \R^3$ with weights all equal to 1.  The curve is $\max\{3x+4,x-y-z,y+z+3\} \cdot X$, the weights are given by the labels.}\label{figure_divisor_example}
\end{figure}

\begin{ex}
 The easiest example of a rational function is a \emph{tropical polynomial} 
$$\varphi(x) = \max\{\gnrt{v_i,x} + \alpha_i; i = 1,\dots,r\}$$
with $v_i \in \Z^n,\alpha_i \in \R$. To this function, we can associate its \emph{Newton polytope}
$$P_\varphi = \conv\{(v_i,\alpha_i); i = 1,\dots,r\} \subseteq \R^{n+1}$$
Denote by $N_\varphi$ its normal fan and define $N_\varphi^1 := N_\varphi \cap \{x: x_{n+1} = 1\}$. Then $N_\varphi^1$ can be considered as a complete polyhedral complex in $\R^n$ and it is easy to see that $\varphi$ is affine linear on each cone of this complex. In fact, each cone in the normal fan consists of those vectors maximizing a certain subset of the linear functions $\gnrt{v_i,(x_1,\dots,x_n)} + \alpha_i \cdot x_{n+1}$ at the same time.

So for any tropical polynomial $\varphi$ and any tropical variety $X$ we can compute an appropriate polyhedral structure on $X$ by intersecting it with $N_\varphi^1$. An example is given in Figure \ref{figure_divisor_example}.
\end{ex}

\begin{algorithm}[ht]
  \caption*{\textbf{polymake example: Computing a divisor}.\\ This computes the divisors displayed in figure \ref{figure_divisor_example}.}
  \flushleft
  \texttt{
  \begin{tabular}{l l}
   \atcommand &\$f = new MinMaxFunction(\\
	      &INPUT\_STRING=\>"max(1,x,y,z,-x,-y,-z)");\\   
   \atcommand &\$x = divisor(linear\_nspace(3),\$f);\\
   \atcommand &\$g = new MinMaxFunction(\\
	      &INPUT\_STRING=\>"max(3x+4,x-y-z,y+z+3)");\\
   \atcommand &\$c = divisor(\$x,\$g);\\
  \end{tabular}}
\end{algorithm}

\subsection{Irreducibility of tropical cycles}

\begin{figure}[ht]
 \centering
  \begin{tikzpicture}
    \matrix[row sep=2mm,column sep=10mm] {
      \draw (0,0) -- (-1,0);
      \draw (0,0) -- (0,-1);
      \draw (0,0) -- (1,1); &
      \draw (-1,0) -- (1,0);
      \draw (0,-1) -- (0,1);
      \draw (-1,-1) -- (1,1);
      \\
    };
  \end{tikzpicture}
  \caption{The curve on the left is irreducible. The curve on the right is reducible and there are several different ways to decompose it.}\label{irred_ex_bad}
\end{figure}
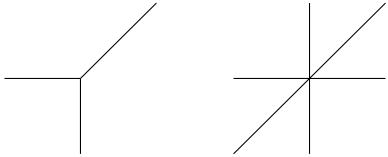

A property of classical varieties that one is often interested in is irreducibility and a decomposition into irreducible components. While one can easily define a concept of irreducible tropical cycles, there is in general no unique decomposition (see Figure \ref{irred_ex_bad}). We can, however, still ask whether a cycle is irreducible and what the possible decompositions are.

\begin{defn}
 We call a $d$-dimensional tropical cycle $X$ \emph{irreducible} if any other $d$-dimensional cycle $Y$ with $\abs{Y} \subseteq \abs{X}$ is an integer multiple of $X$.
\end{defn}

To compute whether a cycle is irreducible, we have to introduce a few notations:

\begin{defn}
 Let $X$ be a tropical cycle with a fixed polyhedral structure $\curly{X}$. Let $N$ be the number of maximal cells $\sigma_1,\dots,\sigma_N$ of $\curly{X}$. We identify an integer vector $\omega \in \Z^N$ with the weight function $\sigma_i \mapsto \omega_i$. We define
 \begin{itemize}
  \item $\Lambda_{\curly{X}} := \{\omega \in \Z^N: (\curly{X},\omega) \textnormal{ is balanced}\}$ (which is a lattice).
  \item $V_{\curly{X}} := \Lambda_\curly{X} \tsr \R$
 \end{itemize}
Now fix a codimension one cell $\tau$ in $\curly{X}$. Let $\curly{S}$ be the induced polyhedral structure of $\Star_\curly{X}(\tau)$. For an integer vector $\omega \in \Z^n$, we denote by $\omega_\curly{S}$ the induced weight function on $\curly{S}$. Then
\begin{itemize}
 \item $\Lambda_\curly{X}^\tau := \{\omega \in \Z^{N}: (\curly{S},\omega_\curly{S}) \textnormal{ is balanced}\}$
 \item $V_\curly{X}^\tau := \Lambda_\curly{X}^\tau \tsr \R$
\end{itemize}
\end{defn}

\begin{remark}
 We obviously have $\Lambda_{\curly{X}} = \bigcap_{\tau \in \curly{X}^{(\dim X -1)}} \Lambda_\curly{X}^\tau$ and similarly for $V_\curly{X}$. Clearly, if $X$ is irreducible, then $\dim V_\curly{X}$ should be 1 and vice versa (assuming that $\gcd(\omega_1,\dots,\omega_N) = 1$, where the $\omega_i$ are the weights on $\curly{X}$). However, so far this definition is tied to the explicit choice of the polyhedral structure. We would like to get rid of this restriction, which we can do using Lemma \ref{irred_lemma_polystructure}. Hence we will also write $V_X$ and $\Lambda_X$. We call $V_X$ the \emphex{weight space} and $\Lambda_X$ the \emphex{weight lattice} of $X$.
\end{remark}

\begin{defn}\label{canonical_defn_equivalence}
 Let $(X,\omega)$ be a $d$-dimensional tropical cycle and $\mathcal{X}$ a polyhedral structure on $X$. We define an equivalence relation on the maximal cells of $\mathcal{X}$ in the following way: Two maximal cells $\sigma, \sigma'$ are equivalent if and only if there exists a sequence of maximal cells $\sigma = \sigma_0,\dots,\sigma_r = \sigma'$, $\sigma_i \in \mathcal{X}^{(d)}$ such that for all $i = 0,\dots,r-1$, the intersection $\sigma_i \cap \sigma_{i+1}$ is a codimension one cell of $\mathcal{X}$, whose only adjacent maximal cells are $\sigma_i$ and $\sigma_{i+1}$.
\end{defn}

\begin{lemma}\label{canonical_lemma_equal}
 Let $(X,\omega)$ be a tropical cycle with polyhedral structure $\curly{X}$ and assume $\sigma,\sigma'$ are equivalent maximal cells of $\curly{X}$. Then: 
\begin{enumerate}
 \item $\omega(\sigma) = \omega(\sigma')$.
 \item If $\omega(\sigma) \neq 0$, then $V_\sigma = V_{\sigma'}$.
\end{enumerate}
\begin{proof}
  \newl
We can assume that $\sigma \cap \sigma' =: \tau \in \mathcal{X}^{(\dim X - 1)}$. 
 \begin{enumerate}
  \item $X$ is balanced at $\tau$ if and only if $\Star_\curly{X}(\tau)$ is balanced, which is a one-dimensional fan with exactly two rays. Such a fan can only be balanced if the weights of the two rays are equal.
  \item Choose any representatives $v_{\sigma/\tau}, v_{\sigma'/\tau}$ of the lattice normal vectors. Then $$\omega(\sigma) v_{\sigma/\tau} + \omega(\sigma') v_{\sigma'/\tau} \in V_\tau$$
Let $g_1,\dots,g_r \in \Lambda^\vee$ such that $$V_\sigma = \ker\begin{pmatrix}g_1\\\vdots\\g_r\end{pmatrix}$$ Since $V_\tau \subseteq V_\sigma$, we have for all $i$:
\begin{align*}
  0 &= g_i(\omega(\sigma) v_{\sigma/\tau} + \omega(\sigma') v_{\sigma'/\tau'}) \\
    &= \omega(\sigma')g_i(v_{\sigma'/\tau'})
\end{align*}
Now $\omega(\sigma') = \omega(\sigma) \neq 0$ implies $v_{\sigma'/\tau} \in V_\sigma$ and since $V_{\sigma'} = V_\tau \times \gnrt{v_{\sigma'/\tau}}$, we have $V_{\sigma'} \subseteq V_\sigma$. The other inclusion follows analogously.
  
 \end{enumerate}
\end{proof}
\end{lemma}

\begin{lemma}\label{irred_lemma_polystructure}
  Let $\curly{X}$ and $\curly{X}'$ be two polyhedral structures of a tropical cycle $X$. Then $V_\curly{X} \cong V_{\curly{X}'}$ (and similarly for $\Lambda_\curly{X} = V_\curly{X} \cap \Z^N$). 
  \begin{proof}
   We can assume without loss of generality that $\curly{X}'$ is a refinement of $\curly{X}$. Denote by $N$ and $N'$ the number of maximal cells of $\curly{X}$ and $\curly{X}'$, respectively and fix an order on the maximal cells of both structures.
First of all, assume two maximal cones of $\curly{X}'$ are contained in the same maximal cone of $\curly{X}$. Since subdividing a polyhedral cell produces equivalent cells in terms of definition \ref{canonical_defn_equivalence}, they must have the same weight by Lemma \ref{canonical_lemma_equal}. Thus the following map is well-defined: We partition $\{1,\dots,N'\}$ into sets $S_1,\dots,S_N$ such that $j \in S_i \iff \sigma_j' \subseteq \sigma_i$ (where $\sigma_j'$ and $\sigma_i$ are maximal cells of $\curly{X}'$ and $\curly{X}$, respectively). Pick representatives $\{j_1,\dots,j_N\}$ from each partitioning set $S_i$ and let $p: V_{\curly{X}'} \to \R^N$ be the projection on these coordinates $j_k$. By the previous considerations, the map does not depend on the choice of representatives. We claim that $\Im(p) \subseteq V_\curly{X}$: Let $\tau$ be a codimension one cell of $\curly{X}$ and $\tau'$ any codimension one cell of $\curly{X}'$ contained in $\tau$. Then $\Star_\curly{X}(\tau) = \Star_{\curly{X}'}(\tau')$, so if $\omega \in \Z^{N'}$ makes $\curly{X}'$ balanced around $\tau'$, then $p(\omega)$ makes $\curly{X}$ balanced around $\tau$. Bijectivity of $p$ is obvious, so $V_\curly{X} \cong V_{\curly{X}'}$.

\end{proof}

\end{lemma}

\begin{theorem}
 Let $(X,\omega)$ be a $d$-dimensional tropical cycle. Then $X$ is irreducible if and only if $g := \gcd(\omega(\sigma), \sigma \in X^{(d)}) = 1$ and $\dim V_X = 1$.
\begin{proof}
  Let $X$ be irreducible. Clearly $g$ must be 1, since otherwise a rational multiple of $X$ would provide a full-dimensional cycle in $X$ not equal to $k \cdot X$ for an integer $k$. Now assume $\dim V_X = \textnormal{rank}(\Lambda_X) > 1$. Then we have an element $\omega' \in \Z^N$ that is not a multiple of $\omega$ and such that $(X,\omega')$ is balanced, which is a contradiction to our assumption that $X$ is irreducible.

  Now let $g = 1 = \dim V_X$. Assume $X$ is not irreducible. Then we can find a polyhedral structure of $X$ and two weight functions $\omega',\omega''$ on this polyhedral structure such that $(X,\omega'), (X,\omega'')$ are both balanced and $\omega' \neq k \cdot \omega''$ for any integer $k$. In particular $\dim V_X \geq 2$, which is a contradiction.
\end{proof}
\end{theorem}

After having laid out these basics, we want to see how we can actually compute this weight space:

\begin{prop}
 Let $\tau$ be a codimension one cell of a $d$-dimensional tropical cycle $X$ in $\R^n$. Let $u_1,\dots, u_k \in \Z^n$ be representatives of the normal vectors $u_{\sigma/\tau}$ for all $\sigma > \tau$. Also, choose a lattice basis $l_1,\dots, l_{d-1}$ of $\Lambda_\tau$. We define the following matrix:
$$ M_\tau := \left( u_1 \;\dots \;u_k \;l_1 \;\dots \;l_{d-1}\right) \in \Z^{n \times (k+d-1)}$$
Then $\Lambda_X^\tau \cong \pi(\ker(M_\tau) \cap \Z^{k+d-1}) \times \Z^{(N-k)}$, where $\pi$ is the projection onto the first $k$ coordinates and $N$ is again the number of maximal cells in $X$.

 \begin{proof}
  Fix an order on the maximal cells of $X$ and let $$J := \{ j \in [N] : \tau \textnormal{ is not a face of } \sigma_j\}.$$ Then clearly $\Z^{(N-k)} \cong \gnrt{e_j; j \in J}_\Z \subseteq \Lambda_X^\tau$ and it is easy to see that $\Lambda_X^\tau$ must be isomorphic to $\Z^{(N-k)} \times \Lambda_{\Star_X(\tau)}$. Hence it suffices to show that $\Lambda_{\Star_X(\tau)}$ is isomorphic to $\pi(\ker(M_\tau) \cap \Z^{k+d-1})$.

  Let $(a_1,\dots,a_k, b_1,\dots,b_l) \in \ker(M_\tau) \cap \Z^{k+d-1}$, Then $\sum a_i u_i = \sum (-b_i) l_i \in \Lambda_\tau$, so $\Star_X(\tau)$ is balanced if we assign weights $a_i$. In particular $(a_1,\dots,a_k) \in \Lambda_{\Star_X(\tau)}$. Since $l_1,\dots,l_{d-1}$ are a lattice basis, any choice of the $a_i$ such that $\Star_X(\tau)$ is balanced fixes the $b_i$ uniquely, so $\pi$ is injective on $\ker(M_\tau)$ and surjective onto $\Lambda_{\Star_X(\tau)}$.
 \end{proof}

\end{prop}


\begin{algorithm}[ht]
  \caption{\sc weightSpace($X$)}
  \begin{algorithmic}[1]
  \INPUT: A pure-dimensional polyhedral complex $X$
  \OUTPUT: Its weight space $V_X$
  \medskip
  \STATE $V_X = \R^N$
  \FOR{$\tau$ a codimension one face of $X$}
    \STATE Compute $M_\tau$ as above
    \STATE $V_X^\tau = \pi(\ker(M_\tau)) + \gnrt{e_j: \tau \textnormal{ is not a face of } \sigma_j}$
    \STATE $V_X = V_X \cap V_X^\tau$
  \ENDFOR
  \RETURN $V_X$
\end{algorithmic}
\end{algorithm}

\begin{algorithm}[ht]
  \caption*{\textbf{polymake example: Checking irreducibility}.\\ This creates the six-valent curve from Figure \ref{irred_ex_bad} and computes its weight space (as row vectors).}
  \flushleft
  \texttt{
  \begin{tabular}{l l}
   \atcommand & \$w = new WeightedComplex(\\
	      &RAYS=\>[[1,0],[1,1],[0,1],[-1,0],[-1,-1],[0,-1]],\\
	      &MAXIMAL\_CONES=\>[[0],[1],[2],[3],[4],[5]],\\
	      &TROPICAL\_WEIGHTS=\>[1,1,1,1,1,1]);\\
   \atcommand & print \$w-\>IS\_IRREDUCIBLE;\\
    0\\
   \atcommand & print \$w-\>WEIGHT\_SPACE;\\
    \multicolumn{2}{l}{1 -1 1 0 0 0}\\
    \multicolumn{2}{l}{0 0 1 0 0 1}\\
    \multicolumn{2}{l}{1 0 0 1 0 0}\\
    \multicolumn{2}{l}{0 1 0 0 1 0}    
  \end{tabular}}
\end{algorithm}

\begin{remark}
 One is often interested in the \emph{positive} weights one can assign to a complex $X$ to make it balanced. This is now very easy using \polymake: Simply intersect $V_X$ with the positive orthant $(\R_{\geq 0})^N$ and you will obtain the \emph{weight cone} of $X$.
\end{remark}

\section{Intersection products in \texorpdfstring{$\R^n$}{Rn}}\label{section_intersection}

There are two main equivalent definitions for a tropical intersection product in $\R^n$, the \emph{fan displacement rule} \cite{fulton}*{Theorem 3.2} and via rational functions \cite{AR}. At first sight, the computationally most feasible one seems to be the latter, since we can already compute it with the means available to us so far:

Let $X,Y$ be tropical cycles in $\R^n$ and $\psi_i = \max\{x_i,y_i\}: \R^n \times \R^n \to \R$. Denote by $\pi: \R^n \times \R^n \to \R^n$ the projection onto the first $n$ coordinates. Then we define $$X \cdot Y := \pi_*(\psi_1 \cdot \dots \cdot \psi_n \cdot (X \times Y))$$
(Here, applying $\pi_*$ just means forgetting the last $n$ coordinates) However, computing this directly turns out to be rather inefficient. The main reason is that, since we compute on the product $X \times Y$, we multiply the number of their maximal cones by each other and double the ambient dimension. As we have discussed earlier, both are factors to which the computation of divisors reacts very sensitively.

A different definition of the intersection product is given by Jensen and Yu: 

\begin{defn}[{\cite{jensenyu}*{Definition 2.4}}]\label{intersection_thm_minkowski}
Let $X,Y$ be tropical cycles in $\R^n$ of dimension $k$ and $l$ respectively. Let $\sigma$ be a $(k+l-n)$-dimensional cone in the complex $X \cap Y$ and $p$ any point in $\relint(\sigma)$. Then $\sigma$ is a cell in $X \cdot Y$ if and only if the Minkowski sum 
$$\Star_X(p) - \Star_Y(p)$$
is complete, i.e.\ its support is $\R^n$. 
\end{defn}

This definition is very close to the fan displacement rule and it is in fact not difficult to see that they are equivalent (\cite{jensenyu}*{Proposition 2.7}). So, at first glance it would seem to be an unlikely candidate for an efficient intersection algorithm. In particular, for $n\geq 6$ it is in general algorithmically undecidable, whether a given fan is complete (see for example the appendix of \cite{fanundecidable}). However, one can also show that $\Star_X(p) - \Star_Y(p)$ can be made into a tropical fan (see \cite{jensenyu}*{Corollary 2.3} for more details). Since $\R^n$ is irreducible, a tropical fan is complete if and only if it is $n$-dimensional. In this case it is a multiple of $\R^n$. 

The weight of the cone $\sigma$ in the above definition is then computed in the following manner:

\begin{defn}[{\cite{jensenyu}}]
 Let $\sigma$ be a polyhedral cell in $X \cdot Y$. Let $p \in \relint(\sigma)$. Then
$$\omega_{X \cdot Y}(\sigma) = \sum_{\substack{\rho_1 \in \Star_X(p), \rho_2 \in \Star_Y(p)\\ \textnormal{s.t. }p \in \relint(\rho_1 - \rho_2)}} \omega_X(\rho_1) \cdot \omega_Y(\rho_2) \cdot ( (\Lambda_{\rho_1} + \Lambda_{\rho_2}) : \Lambda_{\rho_1 + \rho_2})$$
\end{defn}

This now allows us to write down an algorithm based on these ideas:

\begin{algorithm}
 \caption{\sc MinkowskiIntersection}
 \begin{algorithmic}[1]
  \INPUT Two tropical cycles $X,Y$ in $\R^n$ of codimension $k$ and $l$ respectively, such that $k+l \leq n$
  \OUTPUT Their intersection product $X \cdot Y$
  \medskip
  \STATE Compute the $(n-(k+l))$-skeleton $Z$ of $X \cap Y$
  \FOR{$\sigma$ a maximal cell in $Z$}
    \STATE Compute an interior point $p \in \relint \sigma$
    \STATE Compute the local fans $\Star_X(p), \Star_Y(p)$
    \IF{for any $\rho_1 \in \Star_X(p), \rho_2 \in \Star_Y(p)$ the cell $\rho_1 - \rho_2$ is $n$-dimensional}
      \STATE Compute weight $\omega_{X \cdot Y}$ of $\sigma$ as described above
    \ELSE
      \STATE Remove $\sigma$
    \ENDIF
  \ENDFOR
  \RETURN $(Z,\omega_{X \cdot Y})$
 \end{algorithmic}
\end{algorithm}

\begin{algorithm}[ht]
  \caption*{\textbf{polymake example: Computing an intersection product}.\\ This computes the self-intersection of the standard tropical line in $\R^2$.}
  \flushleft
  \texttt{
  \begin{tabular}{l l}
   \atcommand & \$l = tropical\_lnk(2,1);\\
   \atcommand & \$i = intersect(\$l,\$l);\\
   \atcommand & print \$i-\>TROPICAL\_WEIGHTS;\\
   1\\
  \end{tabular}}
\end{algorithm}

\section{Matroid fans}\label{section_matroid}

Matroid fans or \emph{Bergman fans} are an important object of study in tropical geometry, since they are the basic building blocks of what we would consider as \enquote{smooth} varieties. There are several different but equivalent ways of associating a tropical fan to a matroid, see for example \cites{ak, francoisrau, mapoly, troplin, stusolve}. One possibility, which immediately implies a method to compute the fan, is given in \cite{mapoly}*{Proposition 2.5}:

\begin{defn}
 Let $M$ be a matroid on $n$ elements. For $w \in \R^n$ let $M_w$ be the matroid whose bases are the bases $\sigma$ of $M$ of maximal $w$-cost $\sum_{i \in \sigma} w_i$. Then $w$ lies in the \emph{Bergman fan} $B(M)$ if and only if $M_w$ has no loops, i.e.\ the union of its bases is the complete ground set.
\end{defn}

\begin{remark}
 The convex hull of the incidence vectors of the bases of a matroid is a polytope in $\R^n$, the so-called \emph{matroid polytope} $P_M$. So the vectors $w$ maximizing a certain basis are exactly the vectors in the normal cone of the vertex corresponding to that basis. Hence the Bergman fan is a subfan of the normal fan of $P_M$. In addition, we know that it has dimension $\textnormal{rank}(M)$ (this follows immediately from other possible definitions of $B(M)$, see for example \cite{ak}). This gives us an algorithm to compute $B(M)$:
\end{remark}

\begin{algorithm}
 \caption{\sc bergmanFanFromNormalFan}\label{algo_matroid_normal}
 \begin{algorithmic}[1]
  \INPUT A matroid $M$ on $n$ elements, given in terms of its bases.
  \OUTPUT Its Bergman fan in $\R^n$
  \medskip
  \STATE Compute the normal fan $F$ of the matroid polytope $P_M$.
  \STATE $S = $ the $\textnormal{rank}(M)$-skeleton of $F$.
  \FOR{$\xi$ a maximal cone in $S$}
    \STATE Let $\rho$ be the corresponding face of $P_M$ maximized by $\xi$
    \STATE Let $\sigma_1,\dots,\sigma_d$ be the bases corresponding to the vertices of $\rho$.
    \IF{$\bigcup \sigma_i \subsetneq [n]$}
      \STATE Remove $\xi$ from $S$
    \ENDIF
  \ENDFOR
  \RETURN $(S, \omega \equiv 1)$
 \end{algorithmic}
\end{algorithm}

While this algorithm is fairly simple to implement, it is highly inefficient for two reasons: Computing the skeleton of a fan from its maximal cones can be rather expensive, especially if we want to compute a low-dimensional skeleton. But mainly, the problem is that from the potentially many cones of $S$ we often only retain a small fraction. Hence we compute a lot of superfluous information.

\subsection{Computing matroid fans via circuits}

A different definition of a matroid fan can be given in terms of its circuits:

\begin{defn}
 Let $M$ be a matroid on $n$ elements. Then $B(M)$ is the set of all elements $w \in \R^n$ such that for all circuits $C$ of $M$, the minimum $\min\{w_i;\, i \in C\}$ is attained at least twice.
\end{defn}

This definition is used by Felipe Rincón in \cite{tropli} to compute the Bergman fans of linear matroids, i.e.\ matroids associated to matrices. His algorithm requires the computation of a \emph{fundamental circuit} $C(e,I)$ for an independent set $I$ and some element $e \notin I$ such that $I \cup \{e\}$ is dependent:
$$C(e,I) = \{e\} \cup \{i \in I \vert (I \wo \{i\}) \cup \{e\}\textnormal{ is independent}\} $$
It is an advantage of linear matroids that fundamental circuits can be computed very efficiently purely in terms of linear algebra. For general matroids it can still be computed using brute force. With this modified computation of fundamental circuits the algorithm of Rincón can be used to compute Bergman fans of general matroids. It turns out that this is still much faster than the normal fan algorithm above. Table \ref{table_matroid_misc} in the appendix demonstrates this.

\begin{algorithm}[ht]
  \caption*{\textbf{polymake example: Computing matroid fans}.\\ This computes the Bergman fan of a matrix matroid and of the uniform matroid $U_{3,4}$. Both methods use an implementation of the algorithm by Felipe Rincón \cite{tropli}. The first version can only be applied to linear matroids and uses linear algebra to compute fundamental circuits. The second version is for general matroids and significantly slower.}
  \flushleft
  \texttt{
  \begin{tabular}{l l}
   \atcommand & \$m = new Matrix$<$Rational$>$([[1,-1,0,0],[0,0,1,-1]]);\\
   \atcommand & \$bm = bergman\_fan\_linear(\$m);\\
   \atcommand & \$u = matroid::uniform\_matroid(3,4);\\
   \atcommand & \$bm2 = bergman\_fan\_matroid(\$m);
  \end{tabular}}
\end{algorithm}

\subsection{Intersection products on matroid fans}

Intersection products on matroid fans have been studied in \cite{shaw},\cite{francoisrau}. Both approaches however are not suitable for computation. While the approach in \cite{shaw} is more theoretical (except for surfaces, where its approach might lead to a feasible algorithm), the description in \cite{francoisrau} might seem applicable at first. The authors define rational functions which, applied to $B(M) \times B(M)$, cut out the diagonal. Hence they can define an intersection product similar to \cite{AR}. 

However, these rational functions are defined on a very fine polyhedral structure of $B(M)$ induced by its \emph{chains of flats}. These are very hard to compute (\cite{complexitymatroids} gives an \emph{incremental} polynomial time algorithm but \cite{seymour} states that already the number of hyperplanes can be exponential) and the subdivision computed by the algorithm of Rincón is in general much coarser. Also, recall that this approach to computing an intersection product already proved to be inefficient in $\R^n$.

It remains to be seen whether there might be a more suitable criterion for computation of matroid intersection products, maybe similar to Definition \ref{intersection_thm_minkowski}.

\section{Moduli spaces of rational curves}\label{section_moduli}

\subsection{Basic notions}

We only present the basic notations and definitions related to tropical moduli spaces. For more detailed information, see for example \cite{GKM}.

\begin{defn}
 An $n$-\emph{marked rational tropical curve} is a metric tree with $n$ unbounded edges, labelled with numbers $\{1,\dots,n\}$, such that all vertices of the graph are at least trivalent. We can associate to each such curve $C$ its metric vector $(d(C)_{i,j})_{i < j} \in \R^{\binom{n}{2}}$, where $d(C)_{i,j}$ is the distance between the unbounded edges (called \emph{leaves}) marked $i$ and $j$, determined by the metric on $C$.

Define $\Phi_n: \R^n \to \R^{\binom{n}{2}}, a \mapsto (a_i + a_j)_{i < j}$. Then 
$$\mn := \{d(C); C \textnormal{ $n$-marked curve}\} \subseteq \R^{\binom{n}{2}} / \Phi_n(\R^n)$$
is the \emph{moduli space} of $n$-marked rational tropical curves.
\end{defn}

\begin{remark}
 The space $\mn$ is also known as the \emph{space of phylogenetic trees} (\cite{tropgrass}). It is shown (e.g. in \cite{GKM}) that $\mn$ is a pure $(n-3)$-dimensional fan and if we assign weight 1 to each maximal cone, it is balanced (though they do not use the standard lattice, as we will see below). Points in the interior of the same cone correspond to curves with the same \emph{combinatorial type}, i.e.\ forgetting their metric, they are equal. In particular, maximal cones correspond to curves where each vertex is exactly trivalent. We call this particular polyhedral structure on $\mn$ the \emph{combinatorial subdivision}.

The lattice for $\mn$ under the embedding defined above is generated by the rays of the fan. These correspond to curves with exactly one bounded edge. Hence each such curve defines a partition or \emph{split} $I \vert I^c$ on $\{1,\dots,n\}$ and we denote the resulting ray by $v_I$ (note that $v_I = v_{I^c}$). Given any rational $n$-marked curve, each bounded edge $E_i$ of length $\alpha_i$ induces some split $I_i, i=1,\dots d$ on the leaves. In the moduli space, this curve is then contained in the cone spanned by the $v_{I_i}$ and can be written as $\sum \alpha_i v_{I_i}$.

While this description of $\mn$ is very useful to understand the moduli space in terms of combinatorics, it is not very suitable for computational purposes. By dividing out $\Im(\Phi_n)$, we have to make some choice of projection, which would force us to do a lot of tedious (and unnecessary) calculations. Also, the special choice of a lattice would make normal vector computations difficult. However, there is a different representation of $\mn$: It was proven in \cite{ak} and \cite{francoisrau} that
$$\mn \cong B(K_{n-1})/\gnrt{(1,\dots,1)}_\R$$
as a tropical variety, where $K_{n-1}$ is the matroid of the complete graph on $n-1$ vertices. In particular, matroid fans are always defined with respect to the standard lattice. Dividing out the lineality space $\gnrt{(1,\dots,1)}$ of a matroid fan can be done without much difficulty, so we will usually want to represent $\mn$ internally in matroid fan coordinates, while the user should still be able to access the combinatorial information hidden within. 

While the description of $\mn$ as a matroid fan automatically gives us a way to compute it, it turns out that this is rather inefficient. Furthermore, as soon as we want to compute certain subsets of $\mn$, e.g.\ Psi-classes, the computations quickly become infeasible due to the sheer size of the moduli spaces. Hence we would like a method to compute $\mn$ (or parts thereof) in some combinatorial manner. The main instrument for this task is presented in the following subsection.
\end{remark}

\subsection{Prüfer sequences}

Cayley's Theorem states that the number of spanning trees in the complete graph $K_n$ on $n$ vertices is $n^{n-2}$. One possible proof uses so-called \emph{Prüfer sequences}: A Prüfer sequence of length $n-2$ is a sequence $(a_1,\dots,a_{n-2})$ with $a_i \in \{1,\dots,n\}$ (Repetitions allowed!). One can now give two very simple algorithms for converting a spanning tree in $K_n$ into such a Prüfer sequence and vice versa:

Given a spanning tree $T$ in $K_n$, we recursively do the following: Find the smallest leaf $i$ of $T$ and let $v$ be the adjacent node. Attach $v$ to the sequence and remove $i$ from $T$. Continue this until there are only two nodes left.

Given a sequence $P$, let $V = \{1,\dots,n\}$ and recursively do the following: Pick the smallest element $i \in V \wo P$ and let $j$ be the first entry of $P$. Insert an edge between nodes $i$ and $j$. Remove $i$ from $V$ and the first entry from $P$. Continue this until $V$ contains only two elements, then insert an edge between these last two nodes.

%

\begin{figure}[ht]
 \centering
  \begin{tikzpicture}
    \matrix[row sep=2mm,column sep=10mm] {
      \draw (0,0) node[above]{5} -- (-1,1) node[left]{1};
      \draw (0,0) -- (-1,-1) node[left]{3};
      \draw (0,0) -- (1,0) node[above]{6};
      \draw (1,0) -- (2,1) node[right]{2};
      \draw (1,0) -- (2,-1) node[right]{4}; &
      \draw[<->] (0,0) -- (1,0); &
      \draw (0,0) node{(5,6,5,6)};
      \\
    };
  \end{tikzpicture}
  \caption{An example for converting a spanning tree on $K_6$ into a Prüfer sequence and back. The tree can also be considered as a 4-marked rational curve with additional labels at the interior vertices.}\label{pruefer_figure_example}
\end{figure}
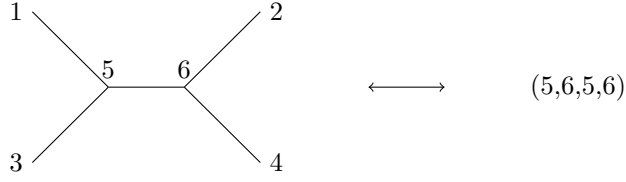

It is easy to see that this induces a bijection (see for example \cite{thebook}*{Chapter 30}). An example for this is given in Figure \ref{pruefer_figure_example}.

As one can see from the picture, tropical rational $n$-marked curves with $d$ bounded edges can also be considered as graphs on $n + d + 1$ vertices: We convert the unbounded leaves into terminal vertices, labelled $1,\dots,n$ and arbitrarily attach labels $n+1,\dots,n+d+1$ to the other vertices. This will allow us to establish a bijection between combinatorial types of rational curves and a certain kind of Prüfer sequence:

\begin{defn}
 A \emph{moduli} Prüfer sequence\index{Prüfer sequence!of rational curves} of order $n$ and length $d$ is a sequence $(a_1,\dots,a_{n+d-1})$ for some $d \geq 0, n \geq 3$ with $a_i \in \{n+1,\dots, n+d+1\}$ such that each entry occurs at least twice.

We call such a sequence \emph{ordered} if after removing all occurrences of an entry but the first, the sequence is sorted ascendingly.

We denote the set of all sequences of order $n$ and length $d$ by $\curly{P}_{n,d}$ and the corresponding ordered sequences by $\curly{P}_{n,d}^<$.
\end{defn}

\begin{ex}
 The sequences $(6,7,8,7,8,6)$ and $(6,7,6,7,8,8)$ are ordered moduli sequences of order 5 and length 2, but the sequence $(6,8,8,7,7,6)$ is not ordered.
\end{ex}

\begin{defn}
 For fixed $n$ and $d$ we call two sequences $p,q \in \curly{P}_{n,d}$ \emph{equivalent} if there exists a permutation $\sigma \in \symm(\{n+1,\dots,n+d+1\})$ such that $q_i = \sigma(p_i)$ for all $i$.
\end{defn}

\begin{remark}
 It is easy to see that for fixed $n$ and $d$ the set $\curly{P}_{n,d}^<$ forms a system of representatives of $\curly{P}_{n,d}$ modulo equivalence, i.e.\ each sequence of order $n$ and length $d$ is equivalent to a unique ordered sequence.

We will need this equivalence relation to solve the following problem: As stated above, we want to associate Prüfer sequences to rational tropical curves by assigning vertex labels $n+1,\dots, n+d+1$ to all interior vertices. There is no canonical way to do this, so we can associate different sequences to the same curve. But two different choices of labellings will then yield two equivalent sequences.
\end{remark}

\begin{prop}\label{pruefer_prop_bijection}
 The set of combinatorial types of $n$-marked rational tropical curves is in bijection to $\bigcup_{d=0}^{n-3} \curly{P}_{n,d}^<$. More precisely, the set of all combinatorial types of curves with $d$ bounded edges is in bijection to $\curly{P}_{n,d}^<$. 
\begin{proof}
 The bijection is constructed as follows: Given an $n$-marked rational curve $C$ with $d$ bounded edges, consider the unbounded leaves as vertices, labelled $\{1,\dots,n\}$. Assign vertex labels $\{n+1,\dots,n+d-1\}$ to the inner vertices. Then compute the Prüfer sequence $P(C)$ of this graph using Prüfer's algorithm and take the unique equivalent ordered sequence as image of $C$.

First of all, we want to see that $P(C) \in \curly{P}_{n,d}$. Since $C$ has $n+d+1$ vertices if considered as above, the associated Prüfer sequence has indeed length $n+d-1$. Furthermore, the $n$ smallest vertex numbers are assigned to the leaves, so they will never occur in the Prüfer sequence. Hence $P(C)$ has only entries in $\left\{n+1,\dots,n+d+1\right\}$. In addition, it is easy to see that each interior vertex $v$ occurs exactly $\val(v) -1$ times (since we remove $\val(v)-1$ adjacent edges before the vertex becomes itself a leaf), i.e.\ at least twice.

Injectivity follows from the fact that if two curves induce the same ordered sequence, they can only differ by a relabelling of the interior vertices, so the combinatorial types are in fact the same. Surjectivity is also clear, since the graph constructed from any $P \in \curly{P}_{n,d}^<$ is obviously a labelling of a rational $n$-marked curve.
\end{proof}
\end{prop}

We now present Algorithm \ref{pruefer_alg_seqToPartitions}, which, given a moduli sequence, computes the corresponding combinatorial type in terms of its edge splits:

\begin{algorithm}
 \caption{\sc combinatorialTypeFromPrueferSequence($P$,$n$)}\label{pruefer_alg_seqToPartitions}
\begin{algorithmic}[1]
  \INPUT A moduli sequence $P = (p_1,\dots,p_N) \in \curly{P}_{n,d}$
  \OUTPUT The rational tropical $n$-marked curve associated to $P$ in terms of the splits $I_1,\dots,I_d$ induced by its bounded edges.
  \medskip
  \STATE $d = N - n +1$
  \STATE $V = \{1,\dots,n+d+1\}$
  \STATE $A_{n+1},\dots,A_{n+d+1} = \emptyset$
  \STATE //First: Connect leaves
  \FOR{$i = 1 \dots n$}
    \STATE $A_{p_i} = A_{p_i} \cup \{i\}$
    \STATE $V = V \setminus \{i\}$
    \STATE $P = (p_{i+1},\dots,p_N)$
  \ENDFOR	
  \STATE //Now create internal edges
  \FOR{$i = n+1 \dots n+d-1$}
    \STATE $v = \min V \setminus P$
    \STATE $I_{i-n} = A_v$
    \IF{\textnormal{length}$(P) > 0$}
      \COMM We denote by $P[0]$ the first element of the sequence $P$.
      \STATE $A_{P[0]} = A_{P[0]} \cup A_v$
      \STATE $V = V \setminus \{i\}$
      \STATE $P = (p_{i+1},\dots,p_N)$
    \ENDIF    
  \ENDFOR
  \STATE //Create final edge
  \STATE $I_d = A_{\min V}$
\end{algorithmic}
\end{algorithm}

\begin{theorem}
 In the notation of Algorithm \ref{pruefer_alg_seqToPartitions} the procedure  generates the set of splits $I_1,\dots,I_d$ of the combinatorial type corresponding to $P$. More precisely: If $v(i)$ is the element chosen from $V$ in iteration $i \in \{n+1,\dots,n+d-1\}$, then $I_{n-i}$ is the split on the leaves $\{1,\dots,n\}$ induced by the edge $\{v(i),p_i\}$.
\begin{proof}
 Let $v(i)$ be the element chosen in iteration $i$, corresponding to vertex $w_i$ in the curve. In particular $i \notin P$. This means that $i$ has already occurred $\val(w_i)-1$ times in the sequence $P$ as $P[0]$. Hence the node $v(i)$ is already $(\val(w_i)-1)$-valent, i.e.\ connected to nodes $q_j; j = 1,\dots,\val(w_i)-1$. If $q_j$ is a leaf (i.e.\ $\leq n$), then $q_j \in A_{v(i)}$. Otherwise, $q_j$ must have been chosen as $v(k)$ in a previous iteration $k < i$. Hence it must already be $\val(w_k)$-valent. Inductively we see that each node ``behind'' $q_j$ is either a leaf or has already full valence. In particular, no further edges will be attached to any of these nodes.

By induction on $i$, the edge $\{v(i),q_j\}$ (assuming $q_j$ is not a leaf) corresponds to the split $A_{q_j}$. In particular, $A_{q_j}$ has been added to $A_v$. Hence $A_v$ is the split induced by the edge $\{v(i),p_i\}$.
\end{proof}
\end{theorem}

\begin{ex}
 Let us apply algorithm \ref{pruefer_alg_seqToPartitions} to the following sequence $P \in \mathcal{P}_{8,5}^<$ (see figure \ref{pruefer_figure_modseq_example} for a picture of the corresponding curve):
$$P = (9,9,10,10,11,11,12,12,13,13,14,14).$$
The algorithm begins by attaching the leaves $\{1,\dots,8\}$ to the appropriate vertices, i.e.\ after the first for-loop we have $A_9 = \{1,2\}, A_{10} = \{3,4\}, A_{11} = \{5,6\},A_{12} = \{7,8\}$ and $P = (13,13,14,14), V = \{9,\dots,14\}$. Now the minimal element of $V \wo P$ is $v = 9$. We set $I_1 = A_9 = \{1,2\}$ to be the split of the first edge. Then we connect the vertex $9$ to the first vertex in $P$, which is 13. Hence $A_{13} = A_{13} \cup A_9 = \{1,2\}$. We remove $9$ from $V$ and set $P$ to be $(13,14,14)$. Now $v = \min V \wo P = 10$. We obtain the second split $I_2 = A_{10} = \{3,4\}$. Then we connect vertex 10 to 13, so $A_{13} = A_{13} \cup A_{10} = \{1,2,3,4\}$. We set $V = \{11,12,13,14\}$ and $P = (14,14)$. In the next two iterations we obtain splits $I_3 = A_{11} = \{5,6\}, I_4 = A_{12} = \{7,8\}$ and we connect both 11 and 12 to 14, setting $A_{14} = \{5,6,7,8\}$. Now $P = ()$ and $V = \{13,14\}$, so we leave the for-loop and set the final split to be $I_5 = A_{13} = \{1,2,3,4\}$.
\begin{figure}[ht]
 \centering
  \begin{tikzpicture}
    \draw (0,0) node[below=2pt]{$9$} -- (-0.5,0.3) node[left]{$1$};
    \draw (0,0) -- (-0.5,-0.3) node[left]{$2$};
    \draw (0,0) -- (1,0) node[below=2pt]{$13$};
    \draw (1,0) -- (1,0.5) node[right]{$10$};
    \draw (1,0.5) -- (0.7,1) node[above]{$3$};
    \draw (1,0.5) -- (1.3,1)node[above]{$4$};
    \draw (1,0) -- (2,0) node[below = 2pt]{$14$};
    \draw (2,0) -- (2,0.5) node[right]{$11$};
    \draw (2,0.5) -- (1.7,1) node[above]{$5$};
    \draw (2,0.5) -- (2.3,1) node[above]{$6$};
    \draw (2,0) -- (3,0) node[below=2pt]{$12$};
    \draw (3,0) -- (3.5,0.3) node[right]{$7$};
    \draw (3,0) -- (3.5,-0.3) node[right]{$8$};
  \end{tikzpicture}
  \caption{The curve corresponding to the moduli sequence $P$, including labels for the interior vertices.}\label{pruefer_figure_modseq_example}
\end{figure}
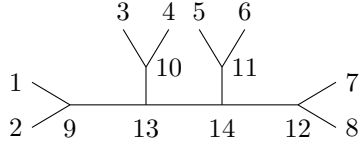
\end{ex}

\subsection{Enumerating  maximal cones of \texorpdfstring{$\mn$}{M0,n}}

We now want to apply the results of the previous section to compute $\mn$. For this it is of course sufficient to compute all maximal cones. More precisely, we will only need to compute all combinatorial types corresponding to maximal cones, i.e.\ rational $n$-marked tropical curves whose vertices are all trivalent. Using Algorithm \ref{pruefer_alg_seqToPartitions} we can then compute its rays $v_{I_1},\dots,v_{I_{n-3}}$. These can easily be converted into matroid coordinates with the construction given in \cite{francoisrau}*{Example 7.2}.

Proposition \ref{pruefer_prop_bijection} directly implies the following:

\begin{corollary}\label{cor_moduli_maximal}
 The maximal cones of $\mn$ are in bijection to all ordered Prüfer sequences of order $n$ and length $n-3$, i.e.\ sequences $(a_1,\dots, a_{2n-4})$ with $a_i$ in $\{n+1,\dots,2n-2\}$ such that each entry occurs exactly twice.
\end{corollary}

This also gives us an easy way to compute the number of maximal cones of $\mn$, which is the \emph{Schröder number}:

\begin{lemma} The number of maximal cones in the combinatorial subdivision of $\mn$ is 
$$(2n-5)!! = \prod_{i=0}^{n-4} (2(n-i) - 5) $$
\begin{proof}
 We prove this by constructing ordered Prüfer sequences of order $n$ and length $n-3$. The sequence has $2n-4$ entries. Since it is ordered, the first entry must always be $n+1$. This entry must occur once more, so we have $2n-5$ possibilities to place it in the sequence. Assume we have placed all entries $n+1, \dots, n+k$, each of them twice. Then the first free entry must be $n+k+1$, since the sequence is ordered and we have $2(n-k)-5$ possibilities to place the remaining one. This implies the formula.
\end{proof}
\end{lemma}

As one can see, the complexity of this number is in $\mathcal{O}(n^{n-3})$, so there is no hope for a fast algorithm to compute all of $\mn$ for larger $n$ (except using symmetries). As we will see later, however, we are sometimes only interested in certain subsets or local parts of $\mn$.

\begin{algorithm}[ht]
  \caption*{\textbf{polymake example: Computing $\mn$}.\\ This computes tropical $\mathcal{M}_{0,8}$ and displays the number of its maximal cones. }
  \flushleft
  \texttt{
  \begin{tabular}{l l}
   \atcommand & \$m = tropical\_m0n(8);\\
   \atcommand & print \$m-\>MAXIMAL\_CONES-\>rows();\\
   10395
  \end{tabular}}
\end{algorithm}

\subsection{Computing products of Psi-classes}

In complex algebraic geometry, Psi-classes on the moduli spaces $\bar{\curly{M}}_{g,n}$ are the first Chern classes of the cotangent bundles of the sections of the universal family (see for example \cite{psinotes} for more details). They became especially interesting, when Witten discovered their relation to string theory and quantum gravity (\cite{witten_psi}). In enumerative geometry, they are useful to count curves satisfying certain tangency conditions. Combining these with pullbacks of evaluation maps to enforce incidence conditions one obtains the so-called descendant Gromov-Witten invariants.

 For the genus 0 case, Mikhalkin suggested a tropical analogon (\cite{mikh_psi}): He defined the $i$-th Psi-class $\psi_i$ as the (closure of the) locus of curves in $\mn$ with a unique four-valent vertex to which the $i$-th leaf is attached. This is a subfan of the codimension one skeleton of $\mn$ and assigning weight 1 to each cone gives a tropical variety. A more detailed study of tropical Psi-classes on $\mn$ was then undertaken in \cite{psiclasses}: The authors describe them as (multiples of) certain divisors of rational functions, but also in combinatorial terms. For nonnegative integers $k_1,\dots,k_n$ and $I \subseteq [n]$, they define $K(I) := \sum_{i \in I} k_i$. Then one of their main results is the following theorem:

\begin{theorem}[{\cite{psiclasses}*{Theorem 4.1}}]
 The intersection product $\psi_1^{k_1}\cdot \dots \cdot \psi_n^{k_n} \cdot \mn$ is the subfan of $\mn$ consisting of the closure of the cones of dimension $n-3 - K([n])$ corresponding to the abstract tropical curves $C$ such that for each vertex $V$ of $C$ we have $\val(V) = K(I_V) + 3$, where
$$I_V = \{i \in [n] : \textnormal{leaf $i$ is adjacent to $V$ } \} \subseteq [n].$$
The weight of the corresponding cone $\sigma(C)$ is
$$\omega(\sigma(C)) = \frac{\prod_{V \in V(C)} K(I_V)!}{\prod_{i=1}^n k_i!}.$$
\end{theorem}

In combination with Proposition \ref{pruefer_prop_bijection}, this allows us to compute these products in terms of Prüfer sequences:

\begin{corollary}
 The maximal cones in $\psi_1^{k_1} \cdot \dots \cdot \psi_n^{k_n} \cdot \mn$ are in bijection to the ordered moduli sequences $P \in \curly{P}_{n,n-3-K([n])}^<$ that fulfill the following condition:

Let $d = n-3-K([n])$ and $k_i = 0$ for $i = n+1,\dots, n+d-1$. For any $a \in \left\{n+1,\dots, n+d+1\right\}$ let $j_1,\dots, j_{l(a)} \in \{1,\dots,n+d-1\}$ be the indices such that $P_{j_i} = a$. Then 
$$l(a) = 2+ \sum_{l = 1}^{l(a)} k_{j_l}$$
\begin{proof}
 Recall that any entry $a$ corresponding to a vertex $v_a$ in the curve $C(P)$ occurs exactly $\val(v_a)-1$ times. By the theorem above the valence of a vertex is dictated by the leaves adjacent to it. Furthermore, the leaves adjacent to a vertex $v_a$ can be read off of the first $n$ entries of the sequence: Leaf $i$ is adjacent to $v_a$ if and only if $P_i = a$.

So, given a curve in the Psi-class product, vertex $v_a$ must have valence $3 + K(I_{v_a})$, so it occurs $2 + K(I_{v_a}) = 2 + \sum_{i: P_i = a} k_i$ times. Conversely, given a sequence fulfilling the above condition, we obviously obtain a curve with the required valences.
\end{proof}
\end{corollary}

We now want to give an algorithm that computes all of these Prüfer sequences. As it turns out, this is easier if we require the $k_i$ to be in decreasing order, i.e.\ $k_1 \geq k_2 \geq \dots \geq k_n$. In the general case we will then have to apply a permutation to the $k_i$ before computation and to the result afterwards. The general idea is that we recursively compute all possible placements of each vertex that fulfill the conditions imposed by the $k_i$ (if we place vertex $a$ at leaf $i$ with $k_i > 0$, then it has to occur more often). Due to its length, the algorithm has been split into several parts: {\sc iteratePlacements} goes through all possible entries of the Prüfer sequence recursively. It uses {\sc placements} to compute all possible valid distributions of an entry, given a certain configuration of free spaces in the Prüfer sequence.

\begin{algorithm}[ht]
 \caption{\sc psiProductSequencesOrdered($k_1,\dots,k_n$)}\label{pruefer_alg_psi}
 \begin{algorithmic}[1]
  \INPUT Nonnegative integers $k_1 \geq k_2 \geq \dots \geq k_n$
  \OUTPUT All Prüfer sequences corresponding to maximal cones in $\psi_1^{k_1} \cdot \dots \cdot \psi_n^{k_n} \cdot \mn$
  \medskip
  \STATE $K = \sum k_i$
  \STATE current\_vertex = $n+1$ 
  \STATE current\_sequence = $(0,\dots,0) \in \Z^{2n-4-K}$
  \STATE exponents = $(k_1,\dots,k_n,0,\dots,0) \in \Z^{2n-4-K}$
  \medskip
  \STATE {{\sc iteratePlacements}(current\_vertex, current\_sequence, exponents)}

 \end{algorithmic}
\end{algorithm}

\begin{algorithm}
 \caption*{ {\sc iteratePlacements}(current\_vertex, current\_sequence, exponents)}
  \begin{algorithmic}[1]
    \IF{current\_vertex $> 2n-2-K$}
      \IF{current\_sequence contains no 0's}
	\STATE append current\_sequence to result
      \ENDIF
    \ELSE
      \STATE $f = \{i: \textnormal{current\_sequence}[i] = 0\}$
      \FOR{$P \in $ {\sc placements($\textnormal{exponents}[i], i \in f$)}}
	\STATE $v = $ current\_sequence
	\STATE Place current\_vertex in $v$ at positions indicated by $P$
	\STATE {\sc iteratePlacements}(current\_vertex+1,$v$,exponents)
      \ENDFOR
    \ENDIF
  \end{algorithmic}
\end{algorithm}

\begin{algorithm}
 \caption*{{\sc placements}($k_1,\dots, k_m$)}
 \begin{algorithmic}[1]
    \INPUT Nonnegative integers $k_1 \geq \dots \geq k_m$
    \OUTPUT All subsets $J \subseteq [m]$ such that $\abs{J} = 2 + \sum_{j\in J} k_j$.
    \medskip
    \IF{$\sum_{i=1}^m k_i > m-2$}
      \RETURN empty list of solutions
    \ENDIF
    \STATE Let $J = \emptyset, i = 1$
    \STATE used$[j] = \emptyset$ for $j = 1,\dots,m$
    \WHILE{$i > 0$}
      \IF{$\abs{J} < 2 + \sum_{j \in J} k_j$}
	\STATE Let $l \in [m] \wo J$ be minimal such that $l > \max J$ and $l \notin $ used$[i]$.
	\IF{There is no such $l$}
	  \STATE \sc{stepDown}
	\ELSE
	  \STATE used$[i$] $=$ used$[i] \cup \{l\}$
	  \STATE $i = i+1$
	  \STATE $J = J \cup \{l\}$
	\ENDIF
      \ELSE
	\IF{$\abs{J} = 2+ \sum_{j \in J} k_j$}
	  \STATE Add $J$ to list of solutions
	\ENDIF
	\STATE \sc{stepDown}
      \ENDIF
    \ENDWHILE
    \RETURN list of solutions
 \end{algorithmic}
\end{algorithm}

\begin{algorithm}
 \caption*{\sc stepDown}
 \begin{algorithmic}[1]
  \STATE used$[i] = \emptyset$
  \STATE $J = J \wo \{\max(J)\}$
  \STATE $i = i -1$
 \end{algorithmic}
\end{algorithm}

\begin{proof}(of Algorithm \ref{pruefer_alg_psi})
 First of all we prove that {\sc placements} computes indeed all possible subsets $J \in [m]$ such that $\abs{J} = 2 + \sum_{j \in J} k_j$. So let $J = \{a_1,\dots,a_N\}$ be such a set with $a_1 \leq \dots \leq a_N$. It is easy to see that in each iteration of the while-loop we have $\abs{J} = i -1$. Let $\delta = (2 + \sum_{j \in J} k_j) - \abs{J}$.

One can see by induction on $\delta$ that, starting in any iteration of the while loop, the algorithm will eventually reach an iteration where $i$ is one smaller. This proves termination of {\sc placements}. 

But we can only reach the iteration where $i = 0$ if in the previous iteration we have tried all indices $\{1,\dots,m\}$ as first element of $J$. In particular, there was a previous iteration, where we chose $l = a_1$ as first element of $J$. Now assume we are in the first iteration where $J = \{a_1,\dots,a_s\}, 1 \leq s < N$. Assuming $\delta > 0$, we can again only decrease $i$ if we have tried all valid placements, including $a_{s+1}$. So assume $\delta = 0$. Then $\{a_1,\dots,a_s\}$ is a valid placement, i.e.\ $s = 2 + \sum_{i=1}^s k_{a_i}$. If we subtract this from the equation for $J$, we obtain
$$0 < N-s = \sum_{i=s+1}^N k_{a_i}$$
In particular, since the $k_i$ are ordered, we must have $k_{a_{s+1}} \geq 1$ and hence also $k_{a_j} \geq 1$ for all $j \leq s$. This implies $$s = 2 + \sum_{i=1}^s k_{a_i} \geq 2 + s$$
which is obviously a contradiction.

With this it is now easy to see that {\sc psiProductSequencesOrdered} computes indeed all the required sequences.
\end{proof}

\begin{ex}
 We do indeed need that $k_1 \geq \dots \geq k_n$ to be able to compute all sequences. Assume $n = 7$ and $(k_1 \dots k_7) = (0,0,0,0,0,1,1)$. A valid sequence would be $(7,7,8,8,9,7,7,9)$, but this sequence would never occur in the algorithm: After having placed the first two 7's in {\sc placements} we would already have $\delta= 0$, so the last two 7's are never tried out.
\end{ex}

For completeness we also give the algorithm for the general case:

\begin{algorithm}
 \caption{{\sc psiProductSequences}($k_1,\dots,k_n$)}
 \begin{algorithmic}[1]
  \INPUT A list of nonnegative integers $\underline{k} = k_1,\dots,k_n$
  \OUTPUT All Prüfer sequences corresponding to maximal cones in $\psi_1^{k_1} \cdot \dots \cdot \psi_n^{k_n} \cdot \mn$
  \medskip
  \STATE Let $\sigma \in \mathcal{S}_n$ such that $\sigma(\underline{k})$ is ordered descendingly
  \STATE $l = $ {\sc psiProductSequencesOrdered($\sigma(\underline{k})$)}
  
  \RETURN $\sigma^{-1}(l)$ (applied elementwise to the first $n$ entries of each sequence)
 \end{algorithmic}

\end{algorithm}

\begin{algorithm}[ht]
  \caption*{\textbf{polymake example: Computing psi classes}.\\ This computes $\psi_1^3 \cdot \psi_2^2 \cdot \psi_6 \cdot \mathcal{M}_{0,9}$ (which is a point) and displays its multiplicity.}
  \flushleft
  \texttt{
  \begin{tabular}{l l}
   \atcommand & \$p = psi\_product(9, new Vector$<$Int$>$(3,2,0,0,0,1,0,0,0));\\
   \atcommand & print \$p-\>TROPICAL\_WEIGHTS;\\
   60\\
  \end{tabular}}
\end{algorithm}

\subsection{Computing rational curves from a given metric}

In previous sections we computed rational curves as elements of the moduli space given by their corresponding bounded edges, i.e.\ the $v_I$ that span the cone containing the curve. Usually, we will be given the curves either in the matroid coordinates of the moduli space or as a vector in $\R^{\binom{n}{2}}$, i.e.\ a metric on the  leaves. It is relatively easy to convert the matroid coordinates to a metric (see \cite{francoisrau}*{Example 7.2}), but it is not trivial to convert the metric to a combinatorial description of the curve, i.e.\ a list of the splits induced by the bounded edges and their lengths.

The paper \cite{buneman} describes an algorithm to obtain a tree from a metric $d$ on a given set $S$ that fulfills the \emph{four-point-condition}, i.e.\ for all $x,y,z,t \in [n]$ we have
$$d(x,y) + d(z,t) \leq \max\{d(x,z) + d(y,t), d(x,t) + d(y,z)\}$$
and \cite{treesproximity}*{Theorem 2.1} shows that the metrics induced by semi-labelled trees (essentially: rational $n$-marked curves) are exactly those which fulfill this condition. 

Note that we can always assume $d(x,y) > 0$ for $x \neq y$ by adding an appropriate element from $\Im(\Phi_n)$. More precisely, if we have an element $d \in \R^{\binom{n}{2}}$ that is equivalent to the metric of a curve modulo $\Im(\Phi_n)$, there is a $k \in \N$ such that $d + k \cdot \Phi_n(\sum e_i)$ is a positive vector fulfilling the four-point-condition. In fact, if $m = d + \sum \alpha_i \Phi_n(e_i)$ is the equivalent metric, then $d + \sum (\alpha_i + \abs{\alpha_i}) \Phi_n(e_i)$ still fulfills the four-point-condition, since adding positive multiples of $\Phi_n(e_i)$ preserves it.

Algorithm \ref{ratcurves_alg_metric} gives a short sketch of the algorithm described in \cite{buneman}*{Theorem 2}. As input, we provide a metric $d$. We then obtain a metric tree with leaves $L$ labelled $\{1,\dots,n\}$ such that the metric induced on $L$ is equal to $d$. This tree corresponds to a rational $n$-marked curve: Just replace the bounded edges attached to the leaf vertices by unbounded edges. It is very easy to modify the algorithm such that it also computes the splits of all edges.

\begin{algorithm}
 \caption{\sc treeFromMetric($d$) \cite{buneman}*{Theorem 2}}\label{ratcurves_alg_metric}
  \begin{algorithmic}[1]
  \INPUT A metric $d$ on the set $[n]$ fulfilling the four-point-condition.
   \OUTPUT A metric tree $T$ with leaf vertices $L$ labelled $\{1,\dots,n\}$ such that the induced metric on $L$ equals $d$.
  \medskip
   \STATE Let $V = \{1,\dots,n\}$
   \WHILE{$\abs{V} > 3$}
    \STATE Find ordered triple of distinct elements $(p,q,r)$ from $V$, such that  $\left.d(p,r) + d(q,r) - d(p,q)\right.$ is maximal
    \STATE Let $t$ be a new vertex and define its distance to the other vertices by
      \begin{align*}d(t,p) &= \frac{1}{2}(d(p,q) - d(p,r) - d(q,r))\\ d(t,x) &= d(x,p) - d(t,p) \textnormal{ for } x \neq p\end{align*}
    \STATE If $d(t,x) = 0$ for any $x$, identify $t$ and $x$, otherwise add $t$ to $V$.
    \STATE Attach $p$ and $q$ to $t$. Then remove $p$ and $q$ from $V$
   \ENDWHILE
   \STATE Compute the tree on the remaining vertices using linear algebra.
  \end{algorithmic}
\end{algorithm}

\begin{algorithm}[ht]
  \caption*{\textbf{polymake example: Converting curve descriptions}.\\ This takes a ray from $\mathcal{M}_{0,6}$ (in its matroid coordinates) and displays it in different representations.}
  \flushleft
  \texttt{
  \begin{tabular}{l l}
   \atcommand & \$m = tropical\_m0n(6);\\
   \atcommand & \$r = \$m-\>RAYS-\>row(0);\\
   \atcommand & print \$r;\\
    \multicolumn{2}{l}{-1 -1 -1 0 -1 -1 0 -1 0}\\
   \atcommand & \$c = rational\_curve\_from\_moduli(\$r);\\
   \atcommand & print \$c;\\
    \multicolumn{2}{l}{(1,2,3,4)} \# This means that this ray represents $v_{\{1,2,3,4\}}$\\
    \atcommand & print \$c-\>metric\_vector;\\
    \multicolumn{2}{l}{0 0 0 1 1 0 0 1 1 0 1 1 1 1 0 \# Read as $d(1,2),d(1,3),\dots,d(5,6)$}\\
   \end{tabular}}
\end{algorithm}
\subsection{Local bases of \texorpdfstring{$\mn$}{M0,n}}

When computing divisors or intersection products on moduli spaces $\mn$, a major problem is the sheer size of the fans, in the number of cones and in the dimension of the ambient space. The number of cones can usually be reduced to an acceptable amount, since one often knows that only a handful of cells is actually relevant. However, the ambient dimension of $\mn$ is $\binom{n}{2} - n = \frac{n^2 - 3n}{2} \in \mathcal{O}(n^2)$. Convex hull computations and operations in linear algebra thus quickly become expensive. We will show, however, that locally at any point $0 \neq p \in \mn$, the span of $\Star_{\mn}(p)$ has a much lower dimension. Hence we can do all our computations locally, where we embed parts of $\mn$ in a lower-dimensional space. Let us make this precise:

\begin{defn}
 Let $\tau$ be a $d$-dimensional cone of $\mn$. We define 
$$V(\tau) := \gnrt{\{\sigma \geq \tau; \sigma \in \mn\}}_\R = \gnrt{U(\tau)}_\R,$$
where $U(\tau) = \bigcup_{\sigma \geq \tau} \relint(\sigma)$. It is easy to see that for any $0 \neq p \in \mn$ and $\tau$ the minimal cone containing $p$, the span of $\Star_{\mn}(p)$ is exactly $V(\tau)$.
\end{defn}

We are now interested in finding a basis for this space $V(\tau)$, preferably without having to do any computations in linear algebra. The idea for this is the following: Let $C_\tau$ be the combinatorial type of an abstract curve represented by an interior point of $\tau$. We want to find a set of rays $v_I$, all contained in some $\sigma \geq \tau$, that generate $V(\tau)$. Each such ray corresponds to separating edges and leaves at a vertex $p$ of $C_\tau$ along a new bounded edge (whose split is of course $I \vert I^c$). We will see that for a fixed vertex $p$ with valence greater than $3$, all the rays separating that vertex span a space that has the same ambient dimension as $\mathcal{M}_{\val(p)}$. In fact, it is easy to see that they must be in bijection to the rays of that moduli space.


Hence the idea for constructing a basis is the following: In addition to the rays of $\tau$, we choose a basis for the \enquote{$\mathcal{M}_{\val(p)}$} at each higher-valent vertex $p$. This choice is similar to the one in \cite{psiclasses}*{Lemma 2.3}. There the authors show that $V_k := \{v_S, \abs{S} = 2, k \notin S\}$ is a generating set of the ambient space of $\mn$ for any $k \in [n]$ and it is easy to see that by removing any element it becomes a basis.

Now fix a vertex $p$ of $C_\tau$ such that $s := \val(p) > 3$. Denote by $I_1,\dots, I_s$ the splits on $[n]$ induced by the edges and leaves adjacent to $p$ (in particular, some of the $I_j$ might only contain one element). We now define
$$W_p := \{v_{I_i \cup I_j}; i,j \neq 1, i \neq j\}$$
(This corresponds to the set $V_1$ described above) and
$$B_p := W_p \wo \{v_{I_2 \cup I_3}\}.$$
Clearly all the following results also hold if we choose $i,j \neq k$ for some $k > 1$ or remove a different element in the definition of $B_p$ (in particular, because the numbering of the $I_i$ is completely arbitrary). To make the proofs more concise, we will however stick to this particular choice. We introduce one final notation: For $\abs{I_j} = 1$, we set $v_{I_j} := 0$.

\begin{lemma}[see also {\cite{psiclasses}*{Lemmas 2.4 and 2.7}}]\label{mnbases_lemma_kmequations}
\newl
\begin{enumerate} 
\item Let $p$ be a vertex of the generic curve $C_\tau$ and define $I_1,\dots,I_s,W_p$ as above. Then
$$\sum_{v \in W_p} v = (s-3) \left(\sum_{j > 1} v_{I_j} \right) + v_{I_1}\equiv 0 \textnormal{ mod } V_\tau.$$
\item Let $v_I$ be a ray in some $\sigma \geq \tau$ and assume it separates some vertex $p$ of $C_\tau$. Define $I_1,\dots,I_s,W_p$ as above. Assume without restriction that $I_1 \subseteq I^c$. Then
$$v_I = \sum_{\substack{v_S \in W_p\\S \subseteq I}} v_S - (m-2)(\sum_{I_j \subseteq I} v_{I_j} + v_I) \equiv \sum_{\substack{v_S \in W_p\\S \subseteq I}} v_S \textnormal{ mod } V_\tau.$$
\end{enumerate}
\begin{proof}\newl
\begin{enumerate}
 \item We define $a = (a_i) \in \R^n$ via $a_i = 1$, if $i \in I_1$ and $a_i = (s-3)$ otherwise. Furthermore we define $b = (b_i) \in \R^n$ via
$$b_i = \begin{cases}
         0, &\textnormal{if } i \textnormal{ is a leaf attached to } p\\
	 1, &\textnormal{if } i \textnormal{ is not a leaf at } p \textnormal{ and lies in } I_1\\
	 (s-3), &\textnormal{if } i \textnormal{ is not a leaf at } p \textnormal{ and does not lie in } I_1.
        \end{cases}$$
We now prove the following equation (to be considered as an equation in $\R^{\binom{n}{2}}$, where each ray is represented by its metric vector):
$$\sum_{v \in W_p} v = (s-3)\left( \sum_{\substack{j > 1 \\ \abs{I_j} > 1}} v_{I_j}\right) + v_{I_1} - \phi_n(b) + \phi_n(a).$$
We index $\R^{\binom{n}{2}}$ by all sets $\mathcal{T} = \{k_1,k_2\}, k_1 \neq k_2$. We have
$$\left( \sum_{v \in W_p} v \right)_{\mathcal{T}} = \begin{cases}
                                                     0, &\textnormal{if } k_1,k_2 \in I_j, j=1,\dots,s\\
						     s-2, &\textnormal{if } k_1 \in I_1,k_2 \in I_j, j > 1 \\
						     2(s-3), &\textnormal{if } k_1 \in I_i, k_2 \in I_j; i,j > 1; i \neq j.
                                                    \end{cases}$$
We now study the right hand side in four different cases:
\begin{enumerate}
 \item If $k_1,k_2 \in I_1$, then both are not leaves at $p$. Hence the right hand side yields $0 + 0 - 2 + 2 = 0$.
 \item If $k_1,k_2 \in I_j, j > 1$, again both are not leaves at $p$. The right hand side now yields $0 + 0 - 2(s-3) + 2(s-3) = 0$.
 \item Assume $k_1 \in I_i, k_2 \in I_j, i,j > 1$ and $i \neq j$. If both are not leaves at $p$, we get $2(s-3) + 0 - 2(s-3) + 2(s-3)$. If only one is a leaf, we get $(s-3) + 0 - (s-3) + 2(s-3)$. Finally, if both are leaves, we get $0 + 0 - 0 + 2(s-3)$. So in any of these cases the right hand side agrees with the left hand side.
 \item Assume $k_1 \in I_1, k_2 \in I_j, j > 1$. If both are not leaves, we get $(s-3) + 1 - (s-3) - 1 + (s-2)$. The other cases are similar.
\end{enumerate}
\item We know that $I$ must be a union of some of the $I_j$ and we assume without restriction that $I = \bigcup_{j \geq k} I_j$ for some $k > 1$. Furthermore we define
$$m := \abs{\{i: I_i \subseteq I\}} = s-k+1.$$
We now prove the following formula (again in $\R^{\binom{n}{2}}$. A similar formula for the representation of a ray $v_I$ in $V_k$ and a similar proof can be found in \cite{psiclasses}*{Lemma 2.7}):
\begin{align}
v_I = & \underbrace{\sum_{i,j \geq k} v_{I_i \cup I_j} - (m-2)\phi_n\left(\sum_{l \in I} e_l\right)}_{=:z} \notag\\
&- (m-2)\underbrace{\left(\sum_{j \geq k} v_{I_j} + v_I - \phi_n\left(\sum_{i=1}^n e_i\right)\right)}_{=:w} \notag\\
 \equiv& \sum_{i,j \geq k} v_{I_i \cup I_j} \textnormal{ mod } V_\tau\label{mnbases_eq_representation}.
\end{align}
To see that the equation holds, let us first compute $w$. We index $\R^{\binom{n}{2}}$ by all sets $\mathcal{T} := \{k_1,k_2\}, k_1 \neq k_2$. Then we have
$$ \left( \sum_{j \geq k} v_{I_j} \right)_{\{k_1,k_2\}} = \begin{cases}
                                                       0, \textnormal{ if } \{k_1,k_2\} \subseteq I_i \textnormal{ for some } i \geq k\\
							1, \textnormal{ if } k_1 \in I, k_2 \notin I \textnormal{ or vice versa}\\
							2, \textnormal{ if } k_1 \in I_i, k_2 \in I_j, i \neq j; i,j \geq k.
                                                      \end{cases}
$$
Hence
\begin{align*}
\left(\sum_{j \geq k} v_{I_j} + v_I \right)_{\{k_1,k_2\}} &= \begin{cases}
                                                 0, \textnormal{ if } \{k_1,k_2\} \subseteq I_i \textnormal{ for some } i \geq k\\
						2, \textnormal{ if } k_1 \in I, k_2 \notin I \textnormal{ or vice versa}\\
						2, \textnormal{ if } k_1 \in I_i, k_2 \in I_j, i \neq j; i,j \geq k
                                                \end{cases}\\
					    &= \begin{cases}
					        0, \textnormal{ if } \{k_1,k_2\} \subseteq I_i \textnormal{ for some } i \geq k\\
						2, \textnormal{ otherwise.}
					       \end{cases}
\end{align*}
Finally we get
$$(w)_{\{k_1,k_2\}} = \begin{cases}
                        -2, \textnormal{ if } \{k_1,k_2\} \subseteq I_i \textnormal{ for some } i \geq k\\
			0, \textnormal{ otherwise.}
                       \end{cases}$$
Thus it remains to prove that
$$\left( v_I - z \right)_{\{k_1,k_2\}} = \begin{cases}
2(m-2),  \textnormal{ if } \{k_1,k_2\} \subseteq I_i \textnormal{ for some } i \geq k\\
0, \textnormal{ otherwise.}                                                                                                                                                                                                                                                             \end{cases}
$$
For this, let $k_1 \neq k_2 \in [n]$. If $\mathcal{T} := \{k_1,k_2\} \subseteq I_i$ for some $i \geq k$, then $(v_I)_{\mathcal{T}} = (v_{I_i \cup I_j})_{\mathcal{T}} = 0$ for all $i,j \geq k$ and $(\phi_n(\sum_{l\in I} e_l))_{\mathcal{T}} = 2$. Thus the formula holds. Now if $k_1 \in I_i, k_2 \in I_j$ for $i \neq j$ and $i,j \geq k$, we still have $(v_I)_\mathcal{T} = 0$. Furthermore, there are $(m-2)$ choices for a ray $v_{I_i \cup I_j'}$  with $j' \neq j$ and $(m-2)$ choices for a ray $v_{I_i' \cup I_j}$ with $i' \neq i$. For these rays, the $\mathcal{T}$-th entry is 1, for all other rays $v_{I_i' \cup I_j'}$ it is 0. Hence $$\left( \sum_{i,j \geq k} v_{I_i \cup I_j}\right)_\mathcal{T} = 2(m-2) = (m-2) \left(\phi_n(\sum_{l\in I} e_l)\right)_{\mathcal{T}}.$$
Finally, if $k_1 \in I$ (say $k_1 \in I_i$), $k_2 \notin I$, then $(v_I)_{\mathcal{T}} = 1$. There are $(m-1)$ choices for a ray $v_{I_i, I_j}$ with $j \neq i$. Since $(\phi_n(\sum_{l\in I} e_l))_{\mathcal{T}} = 1$, we get
$$\left( \sum_{i,j \geq k} v_{I_i \cup I_j}\right)_\mathcal{T} = (m-1) = (m-2) \left(\phi_n(\sum_{l\in I} e_l)\right)_{\mathcal{T}} + 1.$$
Hence equation \ref{mnbases_eq_representation} holds.
\end{enumerate}
\end{proof}
\end{lemma}

\begin{theorem}\label{mnbases_thm_basis}
Let $v_{E_1},\dots,v_{E_t}$ be the rays of $\tau$. Then the set
$$B_\tau := (\bigcup_{\substack{p \in C_\tau^{(0)}\\ \val(p) > 3}} B_p ) \cup \{v_{E_1},\dots,v_{E_t}\}$$
is a basis for $V(\tau)$. In particular, the dimension of $V(\tau)$ can be calculated as
$$\dim V(\tau) = \dim \tau + \sum_{\substack{p \in C_\tau^{(0)}\\ \val(p) > 3}}\left(\binom{\val(p)}{2} - \val(p)\right).$$
\begin{proof}
By Lemma \ref{mnbases_lemma_kmequations} these rays generate $V_\tau$: We can write each $v_I$ in some $\sigma \geq \tau$ in terms of $W_p$ and the bounded edges at the vertex associated to it. The first part of the Lemma then yields that we can replace any occurrence of $v_{I_2 \cup I_3}$ to get a representation in $B_p$ and the bounded edges.

To see that the set is linearly independent, we do an induction on $n$. For $n = 4$ the statement is trivial. For $n > 4$, assume $\tau$ is the vertex of $\mn$. Then $B_\tau$ actually agrees with the set $V_k\wo \{v_S\}$ for some $S$ and we are done. So let $p$ be a vertex of $C_\tau$ that has only one bounded edge attached and denote by $i$ one of the leaves attached to it. It is easy to see that applying the forgetful map $\ft_i$ to $B_\tau$, we get the set $B_{\ft_i(\tau)}$. If $p$ is trivalent, then the ray corresponding to the bounded edge at $p$ is mapped to 0 and all other elements of $B_\tau$ are mapped bijectively onto the elements of $B_{\ft_i(\tau)}$. Since the latter is independent by induction, so is $B_\tau$.

If $p$ is higher-valent, only rays from $B_p$ might be mapped to 0 or to the same element. Hence, if we have a linear relation on the rays in $B_\tau$, we can assume by induction that only the elements in $B_p$ have non-trivial coefficients. But these are linearly independent as well: Let $q$ be any other vertex with only one bounded edge and $j$ any leaf at $q$. $B_p$ is now preserved under the forgetful map $\ft_j$ and hence linearly independent by induction.
\end{proof}
\end{theorem}

At the beginning of this section we introduced the notion that the rays resolving a certain vertex of a combinatorial type $C_\tau$ \enquote{look like $\mathcal{M}_{\val(p)}$}. The results above allow us to make this notion precise:

\begin{corollary}\label{mnbases_cor_decomposition}
 Let $\mathcal{M}$ be any polyhedral structure of $\mn$ (and hence a refinement of the combinatorial subdivision). Let $\tau \in \mathcal{M}$ be a $d$-dimensional cell. Let $C_\tau$ be the combinatorial type of a curve represented by a point in the relative interior of $\tau$. Denote by $p_1,\dots,p_k$ its vertices and by $l$ the number of bounded edges of the curve. Then
$$\Star_{\mn}(\tau) \cong \R^{l-d} \times \mathcal{M}_{\val(p_1)} \times \dots \times \mathcal{M}_{\val(p_k)}$$
\begin{proof}
 First assume $\mathcal{M}$ is the combinatorial subdivision of $\mn$. There is an obvious map $$\psi_\tau: \Star_{\mn}(\tau) \to \mathcal{M}_{\val(p_1)} \times \dots \times \mathcal{M}_{\val(p_k)},$$
defined in the following way: For each vertex $p_i$ of $C_\tau$ fix a numbering of the adjacent edges and leaves, $I_1,\dots,I_{j_i}$. Now for each $v_I$ in some $\sigma \geq \tau$, there is a unique $i \in \{1,\dots,k\}$ such that $v_I$ separates $p_i$. Let $S \subseteq \{1,\dots,j_i\}$ such that $I = \bigcup_{j \in S} I_j$. Again, this choice is unique. Now map $v_I$ to $v_S$ in $\mathcal{M}_{\val(p_i)}$. It is easy to see that this map must be bijective.

First let us see that the map is linear. By Theorem \ref{mnbases_thm_basis} we only have to check that the map respects the relations given in Lemma \ref{mnbases_lemma_kmequations}. But this is clear, since analogous equations hold in $\mn$ (again, see \cite{psiclasses}*{Lemmas 2.4 and 2.7} for details).

For any set of rays $v_{J_1},\dots,v_{J_k}$ associated to the same vertex of $C_\tau$ it is  easy to see that they span a cone in $\mn$ if and only if their images do. Now if $\sigma \geq \tau$ is any cone, we can partition its rays into subsets $S_j, j = 1,\dots,m$ that are associated to the same vertex $p_j$. Each of these sets of rays span a cone $\sigma_j$ which is mapped to a cone in $\mathcal{M}_{\val(p_j)}$. Since $\sigma = \sigma_1 \times \dots \times \sigma_m$, it is mapped to a cone in $\mathcal{M}_{\val(p_1)} \times \dots \times \mathcal{M}_{\val(p_m)}$. Hence $\psi_\tau$ is an isomorphism.

Finally, if $\mathcal{M}$ is any polyhedral structure, let $\tau'$ be the minimal cone of the combinatorial subdivision containing $\tau$. Then $l = \dim \tau'$ and we have
$$\Star_{\mn}(\tau) \cong \R^{l-d} \times \Star_{\mn}(\tau')$$
\end{proof}
\end{corollary}

\begin{algorithm}[ht]
  \caption*{\textbf{polymake example: Local computations in $\mn$}.\\ This computes a local version of $\mathcal{M}_{0,13}$ around a codimension 2 curve $C$ with a single five-valent vertex, i.e.\ it computes all maximal cones containing the cone corresponding to $C$. \atint\ keeps track of the local aspect of this complex, so it will actually consider it as balanced.}
  \flushleft
  \texttt{
  \begin{tabular}{l l}
   \atcommand & \$c = new RationalCurve(N\_LEAVES=\>13,\\
	      & INPUT\_STRING=\>"(2,3) + (2,3,4) + (1,12) + (1,2,3,4,12) + \\
	      & (9,10) + (8,9,10) + (11,13) + (8,9,10,11,13)");\\
   \atcommand & \$m = local\_m0n(\$c);\\
   \atcommand & print \$m-\>MAXIMAL\_CONES-\>rows();\\
   15\\
   \atcommand & print \$m-\>IS\_BALANCED;\\
   1\\
  \end{tabular}}
\end{algorithm}

\newpage
\section{Appendix}\label{section_appendix}

\subsection{Open questions and further research}
 
\subsubsection{More efficient polyhedral computations} As we already discussed in previous sections, most polyhedral operations occurring in tropical computations can in general be arbitrarily bad in terms of performance. However, it remains to be seen how different convex hull algorithms compare in the case of tropical varieties. So far, \atint\ only makes use of the \emph{double-description method} \cite{fukuda_cdd}. Also, measurements indicate that many computations are much faster when all cones involved are simplicial (see the discussion in Section \ref{section_bench_divisor}). Since we are not bound to a fixed polyhedral structure in tropical geometry, it would be interesting to see if we gain anything by subdividing all complexes involved before we do any actual computations.

\subsubsection{A computable intersection product on matroid fans} A very interesting question is whether one can find another description of intersection products in matroid fans that is suitable for computation. One should be able to compute this from the bases of the matroid alone (since that is what is usually given). It would also be interesting to see if there is some local, purely geometric criterion similar to Theorem \ref{intersection_thm_minkowski}.

\subsection{A note on \polymake\ and \atint}

As a tool for polyhedral computations, \polymake\ is available for Linux and Mac under 
\begin{center}
\url{www.polymake.org}
\end{center}Since \atint\ requires a recent version, it is recommended to download the latest package or - even better - source code for installation. One can also find extensive documentation and some tutorials on their webpage.

All the algorithms we discussed in the previous sections have been implemented by the author in \atint, an extension for \polymake. It can be obtained under
\begin{center}
 \url{https://bitbucket.org/hampe/atint}
\end{center}
 Installation instructions and a user manual can be found under the \texttt{Wiki} link. We include here a list of most of the features of the software:
\subsection*{List of features}
\begin{itemize}
\item Creating weighted polyhedral fans/complexes.
\item  Basic operations on weighted polyhedral complexes: Cartesian product, k-skeleton, affine transformations, computing lattice normals, checking balancing condition.
\item Visualization: Display varieties in $\R^2$ or $\R^3$ including (optional) weight labels and coordinate labels (see for example Figures \ref{intro_fig_star} and \ref{figure_divisor_example}).
\item Compute degree and recession fan (experimental).
\item Rational functions: Arbitrary rational functions (given as complexes with function values) and tropical polynomials (min and max).
\item Basic linear arithmetic on rational functions: Compute linear combinations of functions.
\item Divisor computation: Compute (k-fold) divisor of a rational function on a tropical variety.
\item Intersection products: Compute cycle intersections in $\R^n$.
\item Intersection products on matroid fans via rational functions defined by chains of flats (this is very slow).
\item Local computations: Compute divisors/intersections locally around a given face/ a given point.
\item Functions to create tropical linear spaces.
\item Functions to create matroid fans (using a modified version of TropLi by Felipe Rincón \cite{troplisoftware}).
\item Creation functions for the moduli spaces of rational n-marked curves: Globally and locally around a given combinatorial type.
\item Computing with rational curves: Convert metric vectors / moduli space elements back and forth to rational curves, do linear arithmetic on rational curves.
\item Morphisms: Arbitrary morphisms (given as complexes with values) and linear maps
\item Pull-backs: Compute the pull-back of any rational function along any morphism.
\item Evaluation maps. Compute the evaluation map $\textnormal{ev}_i$ on the labelled version of the moduli spaces $\mn$.
\end{itemize}

\subsection{Benchmarks}

All measurements were taken on a standard office PC with 8 GB RAM and 8x2.8 GHz (though no parallelization took place). Time is always given in seconds.

\subsubsection{Divisor computation}\label{section_bench_divisor}

Here we see how the performance of the computation of the divisor of a tropical polynomial on a cycle changes if we change different parameters. In Table \ref{table_divisor_poly} we take a random tropical polynomial $f$ with $l$ terms, where $l \in \{5,10,15\}$. We compute the divisor of this polynomial on $X$, where $X$ is $L \times \R^{k-1}$ in $\R^n$, $L$ the standard tropical line, i.e.\ $X$ has 3 cones. We do this ten times and take the average.

Note that the normal fan of the Newton polytope of $f$ usually has around $l$ maximal cones, but is highly non-simplicial: It can have several hundred rays. If, instead of $f$, we take a polynomial whose Newton polytope is the hypercube in $\R^n$ (here the normal fan is simplicial), then all these computations take less then a second. This is probably due to the fact that convex hull algorithms are very fast on \enquote{nice} polyhedra.

\begin{table}[ht] 
\begin{tabular}{r |r|r|r|r|r|}
  & n=2 &4&6&8&10\\\hline
  k=1 & 0.3 & 0.3 & 0.3 & 0.3 & 0.3 \\
      & 0.2  & 0.3  & 0.5  & 0.4  & 0.5 \\
      & 0.3  & 0.5  & 36.8 & 1242.6 & 2327.3\\\hline
  3   &      & 0.2  & 0.3  & 0.3 & 0.3\\
      &      & 0.4  & 0.6  & 0.5 & 0.5\\
      &      & 0.7  & 38.4 & 1031.7 & 1860.6\\\hline
  5   &      &      & 0.3  & 0.3    & 0.3 \\
      &      &      & 1.9  & 1.8    & 2.2 \\
      &      &      &43.9  & 1519.7 & 2313.5    \\\hline
  7   &      &      &      & 0.4    & 0.4 \\
      &      &      &      & 6.1    & 7.2 \\
      &      &      &      & 2010.9 & 3149.2    \\\hline
  9   &      &      &      &        & 0.4    \\
      &      &      &      &        & 2.4    \\
      &      &      &      &        & 22931    \\
  
\end{tabular}
\vspace{10pt}
\caption{Divisor of a random tropical polynomial with $l = 5,10,15$ terms on $L \times \R^{k-1} \subseteq \R^n$. The time is given in seconds.}\label{table_divisor_poly}
\end{table}

\subsubsection{Intersection products}

We want to see how computation of an intersection product compares to divisor computation. If we apply several rational functions $f_1,\dots,f_k$ to $\R^n$, we can compute $f_1 \cdot \dots \cdot f_k \cdot \R^n$ in two ways: Either as successive divisors $f_1 \cdot (f_2 \cdot (... \cdot \R^n))$ or as an intersection product $(f_1 \cdot \R^n) \cdot \dots \cdot (f_k \cdot \R^n)$. Since successive divisors of rational functions appear in many formulas and constructions, it is interesting to see which method is faster. Table \ref{table_inter_success} compares this for $k = 2$. We take $f$ and $g$ to be random tropical polynomials with 5 terms and average over 50 runs. As we can see, the intersection product is significantly faster in low dimensions, but its computation time grows much more quickly: For $n = 8$, the intersection product takes seven times as long as the divisors.

\begin{table}[ht]
 \begin{tabular}{r| r | r |}
  $n$ & $f \cdot (g \cdot \R^n)$ & $(f \cdot \R^n) \cdot (g \cdot \R^n)$ \\\hline
  3 & 0.62 & 0.14 \\
  4 & 0.68 & 0.24 \\
  5 & 1.04 & 0.38 \\
  6 & 1.42 & 0.84\\
  7 & 1.5  & 2.7 \\
  8 & 1.6  & 11.66\\
 \end{tabular}
\vspace{10pt}
  \caption{Comparing successive divisors to intersection products. Time is given in seconds.}\label{table_inter_success}
\end{table}

\subsubsection{Matroid fan computation}

Here we compare the computation of matroid fans with different algorithms. In Table \ref{table_matroid_moduli}  we compare computation of the moduli space $\mn$, first as the Bergman fan $B(K_{n-1})$ of the complete graph on $n-1$ vertices using the TropLi algorithms, then combinatorially as described in Corollary \ref{cor_moduli_maximal}.

Table \ref{table_matroid_misc} shows some more examples of Bergman fans. We compare the performance of the TropLi algorithms \cite{tropli} to the normal fan Algorithm \ref{algo_matroid_normal}. First, we compute the Bergman fan of the uniform matroid $U_{n,k}$. Note that we compute it as a Bergman fan of a matroid without making use of the matrix structure behind it (the uniform matroid is actually realizable). Then we compute two \emph{linear} matroids, i.e.\ we let TropLi make use of linear algebra to compute fundamental circuits. $C_i$ has as column vectors the vertices of the $i$-dimensional unit cube.
\footnotetext[1]{We did not use the original TropLi program but an implementation of the algorithms in polymake-C++. In the case of linear matroids the original program is actually much faster. This is probably due to the fact that the data types in polymake are larger and that the linear algebra library used by TropLi is more efficient.}


\begin{table}[ht]
 \begin{tabular}{r | r | r |}
  $n$ & TropLi\footnotemark[1] & Cor. \ref{cor_moduli_maximal}\\\hline
  6 & 0 & 0 \\\hline
  7 & 3 & 0 \\\hline
  8 & 921 & 0 \\\hline
 \end{tabular}
\vspace{10pt}
 \caption{Computation of $\mn$}\label{table_matroid_moduli}
\end{table}
%
\begin{table}[ht]
 \begin{tabular}{l | r | r |}
   & TropLi\footnotemark[1] & Algorithm \ref{algo_matroid_normal}\\\hline
  $U_{9,6}$& 0 & 3 \\\hline
  $U_{10,6}$ & 0 & 19 \\\hline
  $U_{11,6}$ & 0 & 87 \\\hline
  $C_3$ & 0 & 1 \\\hline
  $C_4$ & 1 & 29388 \\\hline
 \end{tabular}
\vspace{10pt}
 \caption{Computation of several matroid fans}\label{table_matroid_misc}
\end{table}


\end{document}